\documentclass[11pt]{article}

  \usepackage[letterpaper, left=1in, right=1in, top=1in, bottom=1in]{geometry}
  
  \usepackage[affil-it]{authblk}  % for associating author affiliations
  
  \usepackage{palatino}  % font face
  \usepackage{setspace}  % line spacing
  \usepackage{color}  % font color
  \usepackage[dvipsnames]{xcolor}  % adds more color options

	\usepackage[labelfont=bf]{caption}	% caption style

  \usepackage{amsmath}  % essential math package
  \usepackage{amssymb}  % math symbols
  \usepackage{mathtools}  % more math symbols
  \usepackage{amsthm}  % creating proof and theorem environments
  \allowdisplaybreaks[4]  % breaking multiline equations
  \usepackage{bm}
  \usepackage{bbm}
	\usepackage{relsize}  % can increase the size of math symbols
	\usepackage{algorithm}  % floating environment for algorithmic
	\usepackage{algpseudocode}  % for easy writing of pseudo code
  
  \DeclareMathOperator*{\minimize}{minimize}
  
  \newcommand{\st}{\mathrm{subject\;to}}

	\theoremstyle{definition}

	\newtheorem{example}{Example}
	
	\newtheorem{remark}{Remark}
	
	\usepackage{longtable}  % for tables that span more than one page
	\usepackage{multirow}  % for columns spanning multiple rows
	\usepackage{array}  % for column specifications
	\usepackage{booktabs}  % providing nicer tables, using \toprule, \midrule, etc.
	\usepackage{makecell}
  \usepackage{graphicx}  % enhanced support for graphics, allows e.g. \in­clude­graph­ics
  \usepackage{float}  % better positioning of figures, e.g. by using [h!]
  \usepackage{subfig}  % creating subfigures with \subfloat
  	\graphicspath{{Figures/}}  % path for graphic files

	\usepackage{hyperref}
	\hypersetup{
		pdfstartview = {FitV},  % fits the width of the page to the window
		colorlinks = true,    % false: boxed links; true: colored links
		linkcolor = NavyBlue,  % color of internal links
		citecolor = NavyBlue,  % color of links to bibliography
		filecolor = NavyBlue,  % color of file links
		urlcolor = NavyBlue  % color of external links
	}  

  \usepackage[round]{natbib}  % for natbib, use bibtex as backend
\usepackage{cleveref}
\usepackage{multirow}
\usepackage{amsmath}
\providecommand{\keywords}[1]{\textbf{Keywords: } #1}

\title{Multistage Robust Mixed-Integer Optimization \\ Under Endogenous Uncertainty}
\author[1]{Wei Feng}
\author[1]{Yiping Feng \thanks{Second corresponding author (ypfeng@zju.edu.cn)}}
\author[2]{Qi Zhang \thanks{Primary corresponding author (qizh@umn.edu)}}
\affil[1]{State Key Laboratory of Industrial Control Technology, College of Control Science and Engineering, Zhejiang University, Hangzhou 310027, China}
\affil[2]{Department of Chemical Engineering and Materials Science, University of Minnesota, Minneapolis, MN 55455, USA}
\date{}

\begin{document}

\maketitle

\begin{abstract}
Endogenous, i.e. decision-dependent, uncertainty has received increased interest in the stochastic programming community. In the robust optimization context, however, it has rarely been considered. This work addresses multistage robust mixed-integer optimization with decision-dependent uncertainty sets. The proposed framework allows us to consider both continuous and integer recourse, including recourse decisions that affect the uncertainty set. We derive a tractable reformulation of the problem by leveraging recent advances in the construction of nonlinear decision rules, and introduce discontinuous piecewise linear decision rules for continuous recourse. Computational experiments are performed to gain insights on the impact of endogenous uncertainty, the benefit of discrete recourse, and computational performance. Our results indicate that the level of conservatism in the solution can be significantly reduced if endogenous uncertainty and mixed-integer recourse are properly modeled.
\end{abstract}


\section{Introduction}

Robust optimization has proven to be an effective approach to decision making under uncertainty and has received considerable attention in recent years. While earlier works have only considered the case without recourse \citep{Ben-Tal1998, Ghaoui1998, Bertsimas2004}, lending it a reputation of being overly conservative, effective means have been developed in the last decade to also address the case with recourse, in two- and multistage settings \citep{Ben-Tal2004, Kuhn2011, Zeng2013}. More recent efforts have focused on incorporating not only continuous but also discrete recourse \citep{Bertsimas2007, Hanasusanto2015b, Bertsimas2015a, Postek2016}. For comprehensive reviews of the literature on robust optimization, we refer the reader to \citet{Bertsimas2011}, \citet{Gabrel2014}, and \citet{Yankoglu2019}.

The vast majority of existing works on robust optimization consider \textit{exogenous} uncertainty, which is characterized by fixed uncertainty sets. Significantly fewer have addressed the case of \textit{endogenous}, i.e. decision-dependent, uncertainty although, in practice, many uncertainties are inherently endogenous. Endogenous uncertainty has its origin in stochastic programming \citep{Jonsbraten1998}. Mainly two types of endogenous uncertainty have been considered in the literature: (1) In the case of endogenous uncertainty of \textit{type 1}, decisions affect the realization of the uncertain parameter by altering its underlying probability distribution. For example, a company may shift the probability distribution of the demand for its product toward higher values by lowering the selling price. (2) In the case of endogenous uncertainty of \textit{type 2}, decisions affect the time at which an uncertain parameter materializes or its true value is revealed. The classical example for type-2 endogenous uncertainty is the size of an oilfield for which the true value cannot be determined until one starts drilling and extracting oil from it.

The literature addressing type-1 endogenous uncertainty is relatively sparse. \citet{Ahmed2000} considers single-stage stochastic network problems and uses Luce's choice axiom to develop an expression for the probability of routing along a path, which depends on the network design variables. \citet{Peeta2010} formulate a pre-disaster investment problem in which the failure probabilities of links in a transportation network can be altered by investment decisions related to strengthening those links. Discrete investment decisions are considered, which results in a two-stage stochastic program that allows the choice between a finite number of sets of failure probabilities. \citet{Escudero2018} apply a similar approach to a three-stage resource allocation planning problem for natural disaster relief under type-1 endogenous uncertainty. \citet{Hellemo2018} propose two-stage models with probability distributions that are distorted through an affine transformation or a convex combination of multiple probability distributions. In addition, the authors consider parameterized distributions with the parameters being first-stage decision variables.

A larger number of existing works focus on type-2 endogenous uncertainty. Applications include oil/gas field development \citep{Goel2004}, capacity expansion in process networks \citep{Goel2006}, open-pit mine production scheduling \citep{Boland2008}, clinical trial planning \citep{Colvin2008}, R\&D project portfolio management \citep{Solak2010}, and vehicle routing \citep{Hooshmand2016a}. Most commonly, the problem is formulated as a multistage stochastic program with discrete scenarios and explicit nonanticipativity constraints (NACs) that are active or inactive depending on the decisions related to the endogenous uncertainty. This results in a model that encodes a very large conditional scenario tree, which dramatically increases the computational complexity compared to the case with only exogenous uncertainty. Significant advances have been made in solving such stochastic programs, with approaches that focus on two general strategies: identifying redundant NACs that can be removed \citep{Goel2006, Colvin2008, Gupta2011, Boland2016, Hooshmand2016, Apap2017}, and applying tailored solution algorithms based on Lagrangean decomposition \citep{Goel2006, Gupta2014}, branch-and-cut \citep{Colvin2010}, Benders decomposition \citep{Terrazas-Moreno2012}, sequential scenario decomposition \citep{Apap2017}, or heuristic knapsack decomposition \citep{Christian2015}. \citet{Vayanos2011} apply decision rules to obtain tractable conservative approximations for multistage stochastic programs with type-2 endogenous uncertain parameters that are continuously distributed. We refer to \citet{Apap2017} for a comprehensive review of existing works in this area as well as a unifying framework that addresses problems with both exogenous and type-2 endogenous uncertainty.

Only recently, endogenous uncertainty has also been considered in robust optimization. \citet{Nohadani2018} and \citet{Lappas2018} address static robust optimization, i.e. without recourse, with decision-dependent polyhedral uncertainty sets. Similarly, \citet{Lappas2016} incorporate endogenous uncertainty into a multistage robust process scheduling framework; however, only first-stage decisions can affect the uncertainty set, and all recourse variables are continuous. \citet{Vayanos2019} consider what the authors refer to as decision-dependent information discovery in two- and multistage robust optimization settings using a $K$-adaptability approach; here, the uncertainty set is fixed but the decision maker can decide whether or not to observe some uncertain parameters.

In this work, we consider multistage robust optimization with mixed-integer recourse and decision-dependent uncertainty sets that can be altered in every stage. Specifically, we focus our discussion on polyhedral uncertainty sets that depend linearly on binary decision variables, which only affect the right-hand sides of the inequalities defining the uncertainty set. To derive tractable approximations of the resulting problem, we apply a decision rule approach based on the concept of lifted uncertainty \citep{goh2010distributionally, Georghiou2015a, Bertsimas2018}. This approach has recently been applied to model predictive control \citep{zhang2015convex}, reservoir management \citep{gauvin2017decision}, and transmission expansion planning \citep{dehghan2018multistage1, dehghan2018multistage2}. Here, we further expand the state of the art by introducing discontinuous piecewise linear decision rules for continuous recourse variables. The proposed framework significantly expands our capability to appropriately model endogenous uncertainty in robust optimization settings, with applicability in a variety of areas, such as network design, revenue management, and multiperiod planning. In our computational experiments, we demonstrate the benefits of considering endogenous uncertainty as well as both continuous and binary recourse, and provide results on the computational performance of the proposed models.

The remainder of this paper is organized as follows. In Section \ref{sec:TwoStage}, we present a two-stage robust mixed-integer optimization formulation with endogenous uncertainty, approximate it using decision rules in a lifted space, and derive a tractable mixed-integer linear programming (MILP) reformulation. The proposed approach is then extended to the multistage case in Section \ref{sec:Multistage}. In Section \ref{sec:CaseStudies}, we apply the proposed models to a two-stage design problem and a multistage production planning problem. Finally, in Section \ref{sec:Conclusions}, we close with some concluding remarks.

\paragraph{Notation} We use lowercase and uppercase boldface letters to denote vectors and matrices, respectively, e.g. $\bm{x} \in \mathbb{R}^n$ and $\bm{A} \in \mathbb{R}^{n \times m}$. Scalar quantities are denoted by non-boldface letters. We use $\circ$ and $\mathbbm{1}(\cdot)$ to denote the Hadamard multiplication operator and the indicator function, respectively. Furthermore, $\bm{0}$ and $\bm{e}$ denote the zero and all-ones vectors, respectively, while $\bm{e}_k$ is the standard basis vector whose $k$th element is 1.

\section{The Two-Stage Case}
\label{sec:TwoStage}

From a technical standpoint, the main novelty of this work lies in the incorporation of uncertainty sets that can be affected by binary recourse decisions across multiple stages. This requires a nontrivial integration of decision-dependent uncertainty sets and a means of modeling binary recourse, for which we apply a decision rule approach based on the concept of lifted uncertainty. However, primarily for the sake of clarity, we first discuss the two-stage case in this section, where the uncertainty set only depends on first-stage decisions. Here, we focus on the description of the decision-dependent uncertainty set and the lifting technique that also allows us to introduce discontinuous piecewise linear decisions rules for continuous recourse variables, which are significantly more flexible than traditional affine decision rules. Then, in Section \ref{sec:Multistage}, we turn to the general multistage case, which involves uncertainty sets that depend on binary recourse decisions.

In the following, we present our approach to solving two-stage robust MILPs under endogenous uncertainty that can generally be formulated as follows:
\begin{subequations}
\label{eqn:TwoStage}
\begin{align}
\minimize \quad & x_1 \label{eqn:TwoStageObjective} \\
\st \quad & \bm{\xi}^{\top} \left( \bm{A_n} \bm{x} + \bm{D_n} \bm{y} \right) + \bm{\tilde{a}_n}^{\top} \bm{\tilde{x}}(\bm{\xi}) + \bm{\tilde{d}_n}^{\top} \bm{\tilde{y}}(\bm{\xi}) \leq \bm{\xi}^{\top} \bm{b_n} \quad \forall \, n \in \mathcal{N}, \, \bm{\xi} \in \Xi(\bm{y}) \label{eqn:TwoStageCon} \\
& \bm{x} \in \mathbb{R}^{P}, \; \bm{y} \in \{0,1\}^{Q} \label{eqn:TwoStageDomain1} \\
& \bm{\tilde{x}}(\bm{\xi}) \in \mathbb{R}^{\widetilde{P}}, \; \bm{\tilde{y}}(\bm{\xi}) \in \{0,1\}^{\widetilde{Q}} \quad \forall \, \bm{\xi} \in \Xi(\bm{y}), \label{eqn:TwoStageDomain2}
\end{align}
\end{subequations}
where $\bm{x}$ and $\bm{y}$ are the first-stage variables, and $\bm{\tilde{x}}$ and $\bm{\tilde{y}}$ are the second-stage variables, which are functions of the uncertain parameters $\bm{\xi} \in \mathbb{R}^K$. For ease of exposition, we assume that the objective function \eqref{eqn:TwoStageObjective} is certain, which is without loss of generality since $x_1$ can simply be the auxiliary variable introduced in an epigraph reformulation. The set of constraints is denoted by $\mathcal{N}$, and as stated in \eqref{eqn:TwoStageCon}, they have to hold for all $\bm{\xi}$ in an uncertainty set $\Xi$. Variable domains are specified in constraints \eqref{eqn:TwoStageDomain1} and \eqref{eqn:TwoStageDomain2}. Note that we assume fixed recourse and that all discrete variables are binary.

Formulation \eqref{eqn:TwoStage} implies that the uncertainty set $\Xi$ depends on the binary decisions $\bm{y}$. Indeed, similar to \citet{Nohadani2018}, we consider decision-dependent uncertainty sets of the following form:
\begin{equation}
\label{eqn:UncertaintySet}
\Xi(\bm{y}) = \left\lbrace \bm{\xi} \in \mathbb{R}^K: \bm{W} \bm{\xi} \leq \bm{U} \bm{y}, \; \xi_1 = 1 \right\rbrace,
\end{equation}
which is assumed to be a compact polyhedron for any feasible $\bm{y}$, with $\bm{W} \in \mathbb{R}^{M \times K}$ and $\bm{U} \in \mathbb{R}^{M \times Q}$. We assume that $\xi_1 = 1$, which eases the notation when introducing constant terms in the linear decision rules as shown in Section \ref{sec:DecisionRules}.

\begin{example}
Consider the following example of a decision-dependent uncertainty set:
\begin{equation}
\label{eqn:UncertaintySetExample}
\Xi(\bm{y}) = \left\lbrace
  \; \bm{\xi} \in \mathbb{R}_+^3: 
  \begin{array}{l}
  \xi_2 \leq 7 y_1 + 8 y_2 \\
  \xi_3 \leq 13 y_2 \\
  -8\xi_1 -\xi_2 +2\xi_3 \leq 15y_2 \\
  -13\xi_1 +\xi_2 + \xi_3 \leq 7 y_1 + 2 y_2 \\
  25 \xi_1 + 4 \xi_2 -7 \xi_3 \leq 21 y_1 +  11y_2 \\
  40 \xi_1 - 8 \xi_2 - 3 \xi_3 \leq 0 \\
  \xi_1 = 1
  \end{array}
\right\rbrace.
\end{equation}
Here, $\xi_3$ can only be nonzero if $y_2 = 1$, and $\xi_2$ is forced to be 0 when $y_1$ and $y_2$ are both 0.
% whereas $\xi_3$ is always observed, which means that $z_3$ is fixed to be 1. Hence, $z_3$, just like $z_1$, is omitted from the formulation of the uncertainty set. 
Figure \ref{fig:Example1} shows the projections of the uncertainty set $\Xi$ onto the two-dimensional $(\xi_2,\xi_3)$-space for different values of $y_1$ and $y_2$. One can see that the choice of $\bm{y}$ affects the facets of $\Xi$, leading to different uncertainty sets. In this case, $\Xi$ is largest when $y_1=1$ and $y_2=1$ (see Figure \ref{fig:Example1a}); in fact, $\Xi(1,1)$ is a superset of $\Xi(0,1)$ (Figure \ref{fig:Example1b}) and $\Xi(1,0)$ (Figure \ref{fig:Example1c}), with the latter being merely the line for which $5 \leq \xi_2 \leq 7$. Note that it is not generally true that there always exists a $\Xi(\bm{y}^*)$ such that $\Xi(\bm{y}^*)$ is a superset of all $\Xi(\bm{y})$ for $\bm{y} \neq \bm{y}^*$. Furthermore, in this example, $\Xi$ is an empty set for $y_1=0$ and $y_2=0$; hence, there has to be a constraint in problem \eqref{eqn:TwoStage}, e.g. $y_1 + y_2 \geq 1$, that excludes this solution.

\begin{figure}[ht]\centering
\subfloat[$y_1=1$, $y_2=1$]{
  \label{fig:Example1a}
	\fbox{\includegraphics[width=2.1in]{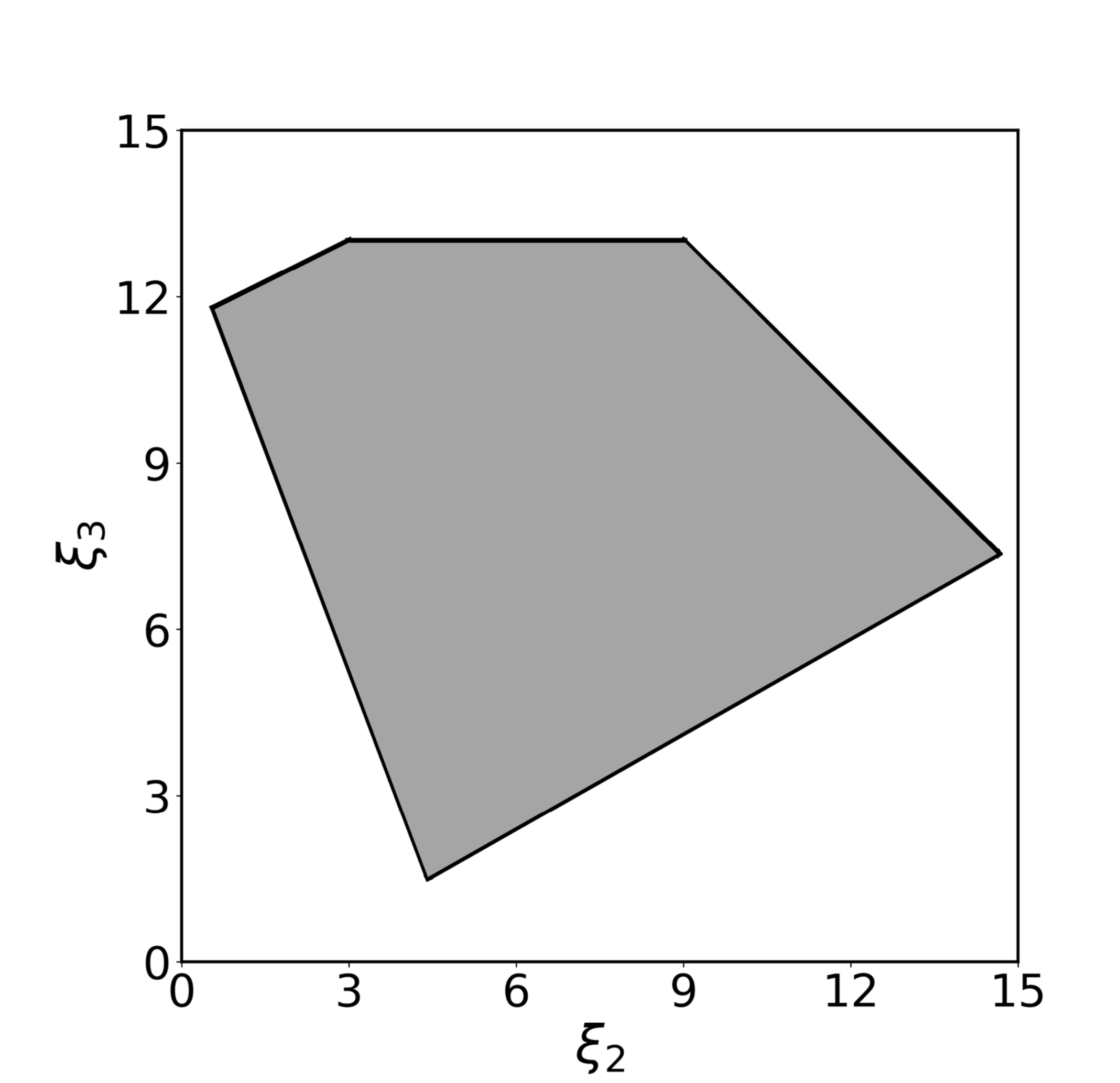}}} 
\subfloat[$y_1=0$, $y_2=1$]{
  \label{fig:Example1b}
	\fbox{\includegraphics[width=2.1in]{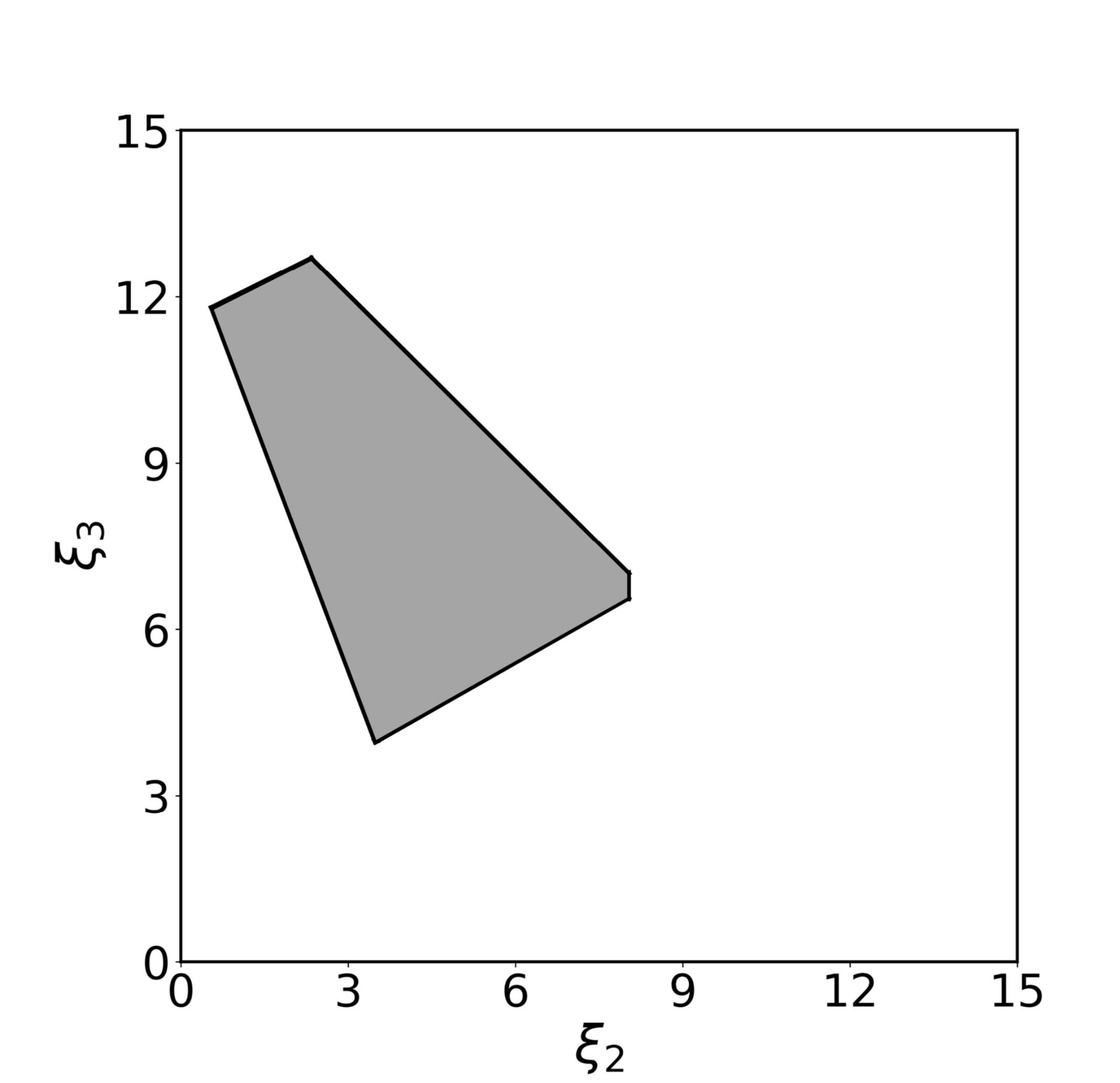}}} 
\subfloat[$y_1=1$, $y_2=0$]{
  \label{fig:Example1c}
	\fbox{\includegraphics[width=2.1in]{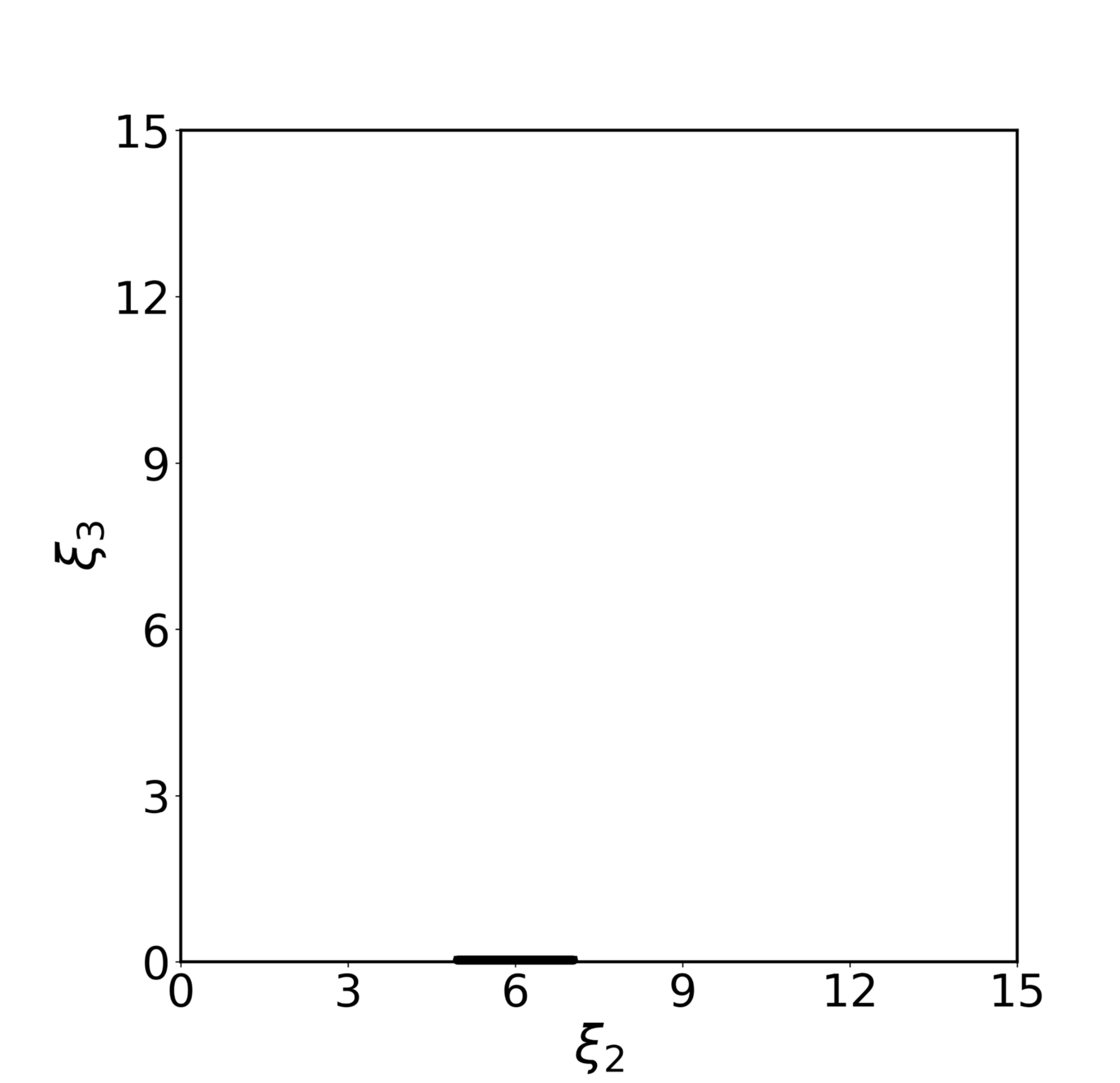}}}
\caption{Projections of $\Xi(\bm{y})$ onto the $(\xi_2,\xi_3)$-space for different values of $y_1$ and $y_2$.}
\label{fig:Example1}
\end{figure}

Clearly, the decision-dependent uncertainty set given by \eqref{eqn:UncertaintySetExample} captures type-1 endogenous uncertainty as the uncertain parameters' support changes with $y$. Depending on the constraints of the problem, it can also consider type-2 endogenous uncertainty. For example, if $y_2 = 0$ renders uncertain parameter $\xi_3$ irrelevant for the problem, e.g. if it is multiplied by $y_2$ everywhere it appears in the model, then the uncertainty set can encode the case in which $\xi_3$ only materializes if $y_2 = 1$. In general, we assume that an uncertain parameter whose materialization is decision-dependent can take the value zero, which is without loss of generality as it can always be achieved with a simple linear translation.

\begin{remark}
In this work, we assume that the uncertainty set can only be affected by binary variables, which allows the final reformulation to be an MILP. The proposed approach can also be applied if $\bm{y}$ are continuous; however, in that case, we would arrive at a nonconvex mixed-integer nonlinear program (MINLP). Furthermore, the approach can be extended to consider random recourse and polyhedral uncertainty sets in which decisions also affect the left-hand sides; however, this would result in a substantial increase in computational complexity and is hence not discussed in this work.
\end{remark}

\end{example}

% It lifts the original uncertainty $\xi_i \in \mathbb{R}$ for $i=2, \ldots, K$ to a $(r_i-1)$-dimensional vector $\bm{\tilde{\xi}}_i$ whose elements are either 0 or 1 as shown by. Due to the degenerate nature of $\xi_1$ in \eqref{eqn:UncertaintySet}, its lifted uncertainty by $\widetilde{L}$ 
% Thus, the dimension of its lifted space is $\widetilde{K} = \sum$

\subsection{Lifted Uncertainty Set}
\label{sec:LiftedUncertainty}

We apply the lifting technique proposed by \citet{Georghiou2015a} to derive tractable decision rules. The basic idea is to lift the original uncertain parameters onto a higher-dimensional space such that linear decision rules can be applied to the lifted parameters, allowing the construction of more flexible and binary decision rules.

\begin{figure}[ht]\centering
\subfloat[$r_i > 1$, $j = 1$]{
  \label{fig: Continuous operator 1}
	\fbox{\includegraphics[height=1in]{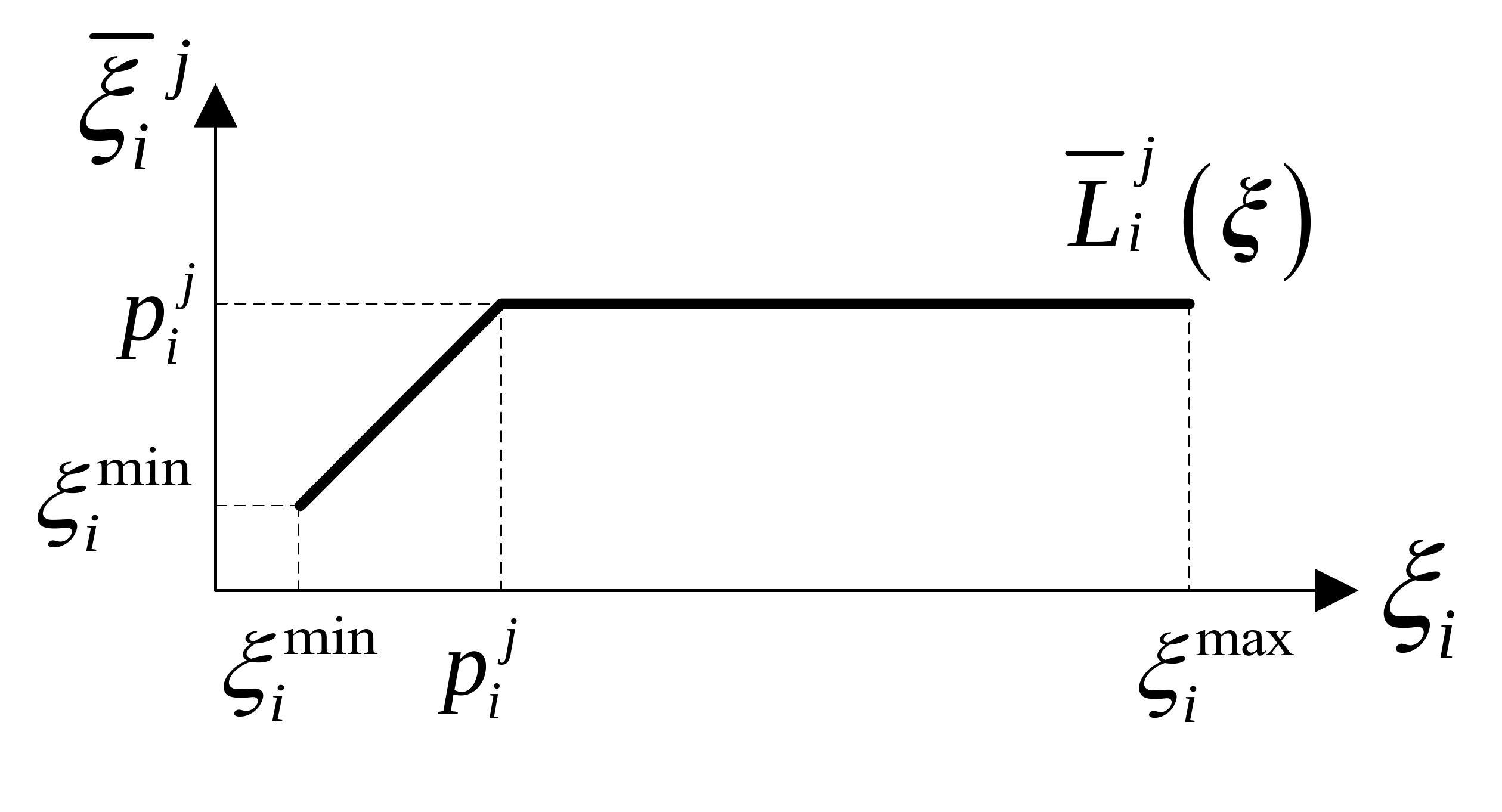}}} 
\subfloat[$r_i > 1$, $j = 2, \ldots, r_i - 1$]{
  \label{fig: Continuous operator 2}
	\fbox{\includegraphics[height=1in]{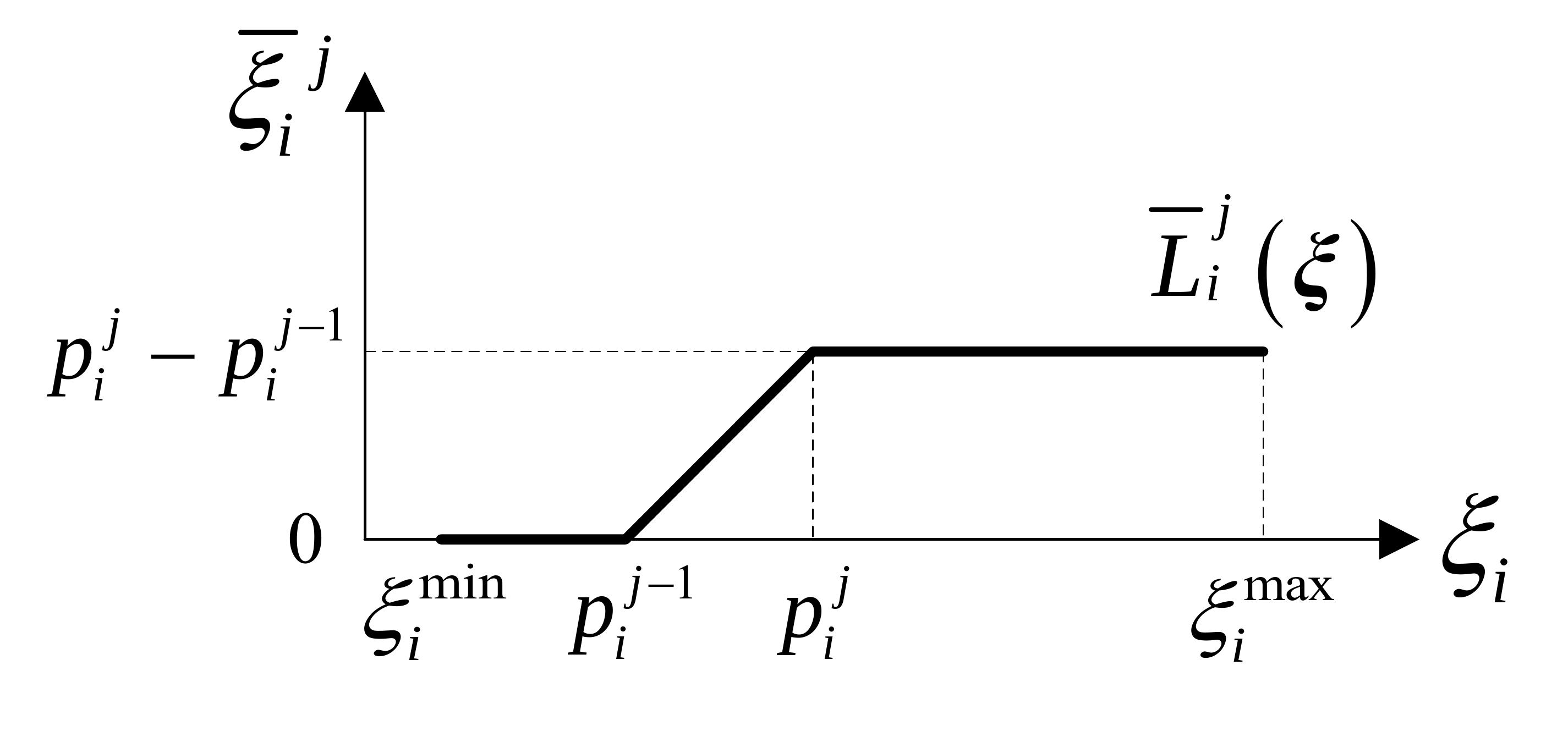}}} 
\subfloat[$r_i > 1$, $j = r_i$]{
  \label{fig: Continuous operator 3}
	\fbox{\includegraphics[height=1in]{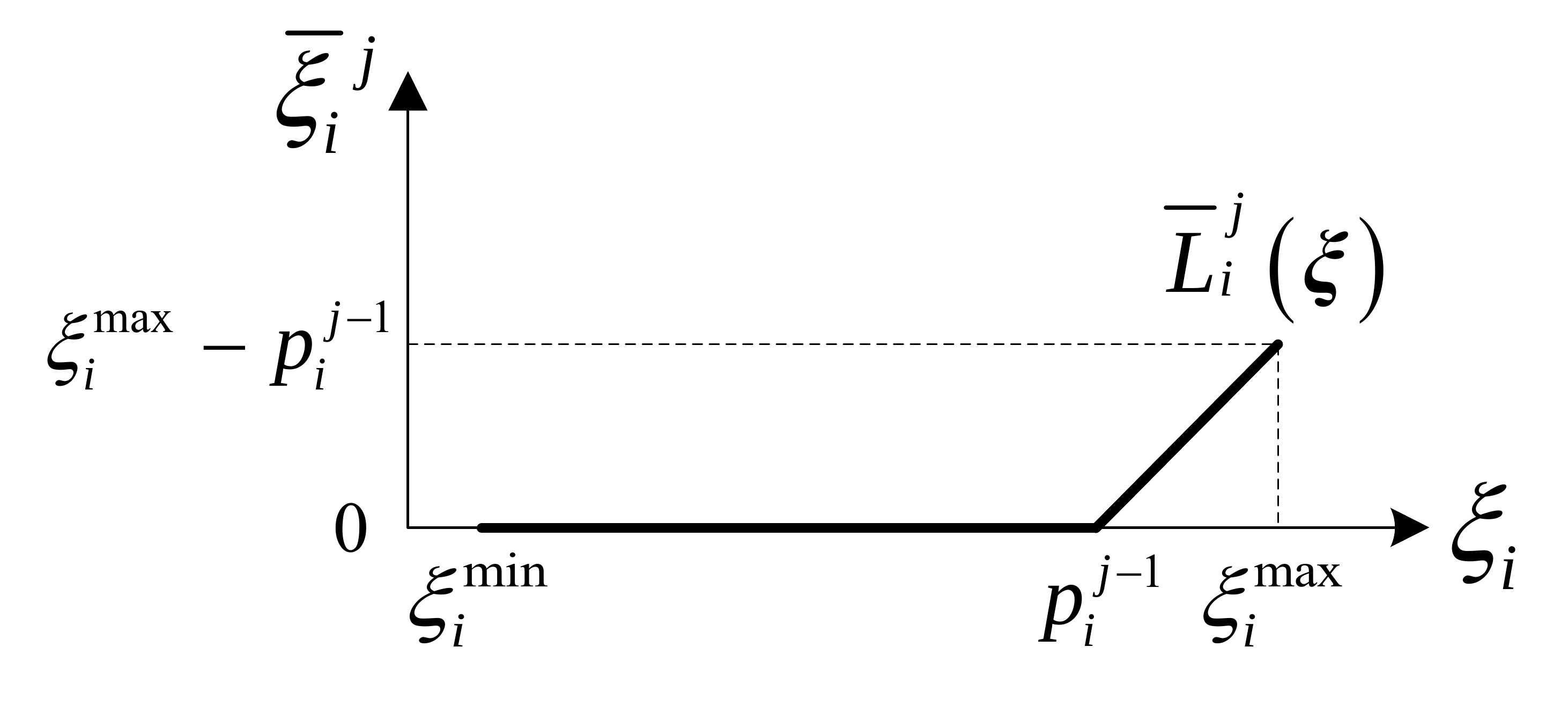}}}
\caption{Illustration of the piecewise linear function $\bar{L}_i^j(\bm{\xi})$.}
\label{fig: Continuous operator}
\end{figure}

We first insert $r_i-1$ breakpoints into the marginal support of each uncertain parameter $\xi_i$ such that
\begin{equation*} \label{eqn: sorting breakpoints}
    \xi_i^{\min} < p_i^1 < p_i^2 < \cdots < p_i^{r_i - 1} < \xi_i^{\max} \quad \forall \, i = 1, \ldots, K,
\end{equation*}
where $\xi_i^{\min}$ and $\xi_i^{\max}$ are the lower and upper bounds of $\xi_i$, respectively. As such, $\{\bm{\xi} \in \mathbb{R}^K: \bm{\xi}^{\min} \leq \bm{\xi} \leq \bm{\xi}^{\max}\}$ is the smallest hyperrectangle that contains $\Xi(\bm{y})$ for all feasible $\bm{y}$. We now define a lifting operator $\bar{L}: \mathbb{R}^{K} \mapsto \mathbb{R}^{\overline{K}}$ that maps the original uncertain parameters $\bm{\xi} \in \mathbb{R}^K$ onto a $\overline{K}$-dimensional space with $\overline{K} =\sum_{i=1}^K r_i$. The vector of lifted uncertain parameters is denoted by $\bm{\bar{\xi}} = \left ( \bm{\bar{\xi}}_1, \ldots, \bm{\bar{\xi}}_K \right ) \in \mathbb{R}^{\overline{K}}$ with $\bm{\bar{\xi}}_i = \left ( \bar{\xi}_i^1, \ldots, \bar{\xi}_i^{r_i} \right ) \in \mathbb{R}^{r_i}$, and the lifting operator $\bar{L} = \left( \bar{L}_1, \ldots, \bar{L}_K \right)$ with $ \bar{L}_i = \left( \bar{L}_i^1, \ldots, \bar{L}_i^{r_i} \right)$ is defined as follows:
\begin{equation} \label{eqn: continuous operator}
    \bar{\xi}_i^j = \bar{L}_i^j(\bm{\xi}):= 
    \left \{
        \begin{array}{lcl}
         {\xi_i}  & \text{if} & r_i = 1 \\
         {\inf \left \{ \xi_i, p_i^j \right\} } & \text{if} & r_i > 1, \, j = 1 \\
         {\sup \left \{ \inf \left \{ \xi_i, p_i^j \right \} - p_i^{j-1}, 0 \right \}} & \text{if} & r_i > 1, \, j = 2, \ldots, r_i - 1 \\
         {\sup \left \{ \xi_i - p_i^{j-1}, 0 \right \}} & \text{if} & r_i > 1, \, j=r_i.
        \end{array}  
    \right .
\end{equation}
As illustrated in Figure \ref{fig: Continuous operator}, $\bar{L}_i^j(\bm{\xi})$ is a piecewise linear function of $\xi_i$. By construction, the original and lifted uncertain parameters have the following relationship:
\begin{equation*} \label{eqn: inverse operator 1}
    \xi_i = \bm{e}^{\top}\bm{\bar{\xi}}_i \quad \forall \, i = 1, \ldots, K.
\end{equation*}
% there exists an inverse operator $\overline{R}: \mathbb{R}^{\overline{K}} \mapsto \mathbb{R}^{K}, \, \bm{\overline{\xi}} \mapsto \bm{\xi}$ such that,
% \begin{equation}
%     \bm{\xi} = \overline{R}(\bar{L}(\bm{\xi})) \quad \forall \bm{\xi} \in \Xi.
% \end{equation}

\begin{figure}[ht]\centering
\includegraphics[width=2in]{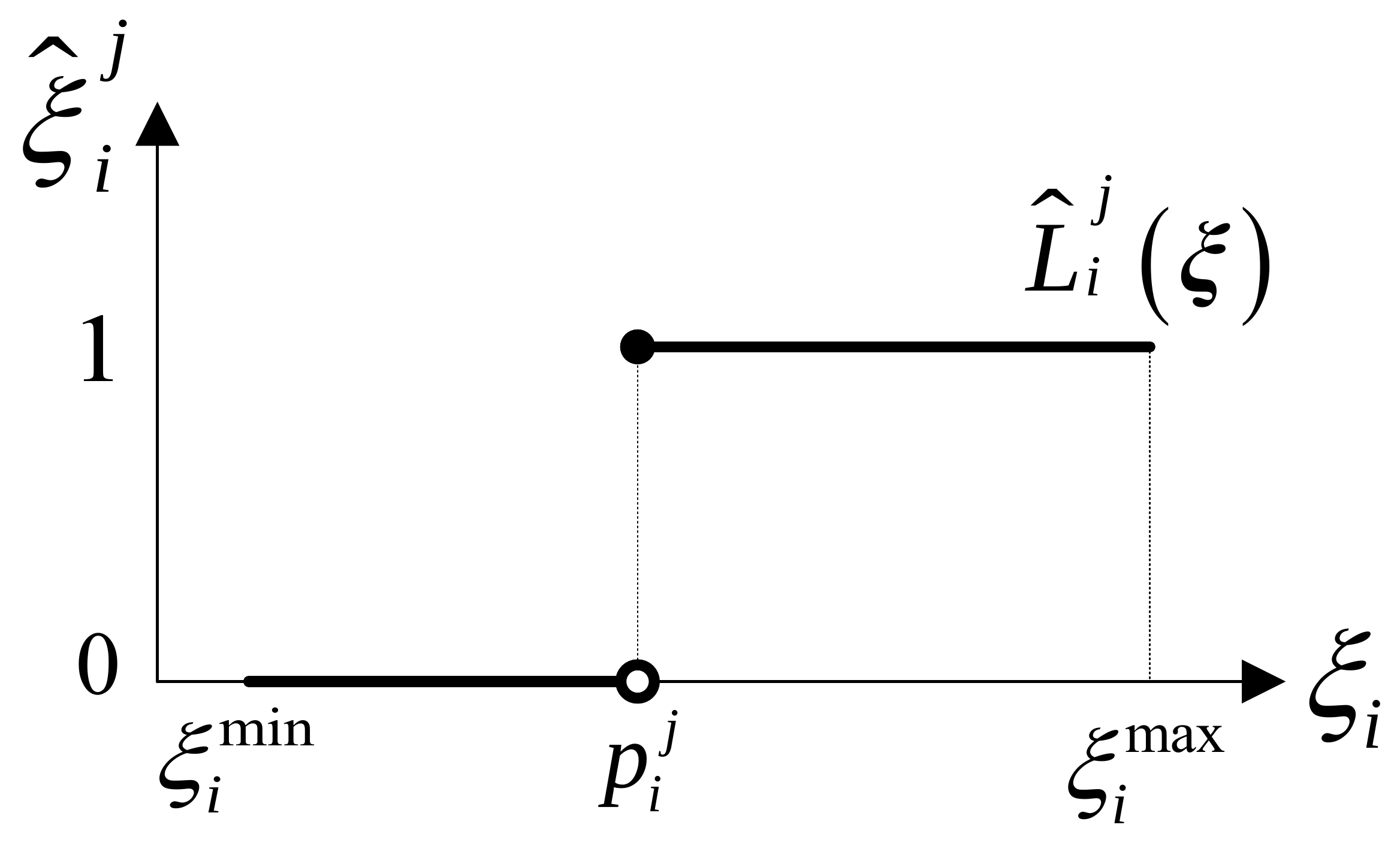}
\caption{Illustration of the piecewise constant function $\widehat{L}_i^j(\bm{\xi})$.}
\label{fig: BinaryOperator}
\end{figure}

Piecewise linear decision rules for continuous recourse variables can be derived using lifted uncertain parameters defined by $\bar{L}$. To also allow binary recourse, we define a second lifting operator $\widehat{L}: \mathbb{R}^{K} \mapsto \{0,1\}^{\widehat{K}}$ as proposed by \citet{Bertsimas2018}. Here, we apply the same breakpoints as in \eqref{eqn: continuous operator} and have $\widehat{K} = \sum_{i=1}^{K} g_i$ with $g_i = \max\{ 1, \, r_i - 1 \}$. We introduce another vector of lifted uncertain parameters $\bm{\hat{\xi}} = \left ( \bm{\hat{\xi}}_1, \ldots, \bm{\hat{\xi}}_K \right ) \in \{ 0, 1 \}^{\widehat{K}}$ with $\bm{\hat{\xi}}_i = \left (\hat{\xi}_i^1, \ldots, \hat{\xi}_i^{g_i} \right) \in \{ 0, 1 \}^{g_i}$, and define $\widehat{L} = \left (\widehat{L}_1, \ldots, \widehat{L}_K \right )$ with $\widehat{L}_i = \left( \widehat{L}_i^1, \ldots, \widehat{L}_i^{g_i} \right)$ as follows:
\begin{equation} \label{eqn: binary operator}
    \hat{\xi}_i^j = \widehat{L}_i^j(\bm{\xi}) := 
    \left \{
        \begin{array}{lcl}
         {1}  & \text{if} & r_i = 1 \\
         {\mathbbm{1}(\xi_i \geq p_i^j)} & \text{if} & r_i > 1, \, j = 1, \ldots, g_i,
        \end{array}  
    \right .
\end{equation}
which is illustrated in Figure \ref{fig: BinaryOperator}.

We introduce a new vector of lifted uncertain parameters, $\bm{\xi'} = \left (\bm{\bar{\xi}}, \, \bm{\hat{\xi}} \right) \in \mathbb{R}^{K'}$ with $K' = \sum_{i=1}^K r_i + g_i$. Note that $\xi_1 = 1$, hence $r_1 = 1$, $\bar{L}_1^1 (\bm{\xi}) = 1$, and $\widehat{L}_1^1 (\bm{\xi}) = 1$. The lifted uncertainty set is then defined as follows:
\begin{equation}
\label{eqn: lifted uncertainty set, open}
    \widetilde{\Xi}'(\bm{y}) =  \left \{ \bm{\tilde{\xi}} = 
    \left( \bm{\xi}, \, \bm{\xi'} \right) 
    = \left( \bm{\xi}, \, \bm{\bar{\xi}}, \, \bm{\hat{\xi}} \right) \in \mathbb{R}^{K + K'}:
    \begin{array}{l}
         \bm{\xi} \in \Xi(\bm{y}) \\
         \bm{\bar{\xi}} = \bar{L}\left(\bm{\xi} \right) \\
         \bm{\hat{\xi}}= \widehat{L} \left( \bm{\xi} \right) 
    \end{array}
    \right \}.
\end{equation}
Due to the discontinuity of $\bm{\hat{\xi}}$ at each breakpoint (see Figure \ref{fig: BinaryOperator}), $\widetilde{\Xi}'$ is an open set. Let $\dot{\widetilde{\Xi}} := \mathrm{cl}(\widetilde{\Xi}')$ be the closure of $\widetilde{\Xi}'$, $\dot{\widetilde{\Xi}}_{\bm{\xi}'}$ denote the projection of $\dot{\widetilde{\Xi}}$ onto the space of $\bm{\xi}'$, and the second-stage variables be functions of the lifted uncertain parameters $\bm{\xi}'$. The two-stage problem then becomes
\begin{subequations}
\label{eqn:TwoStage with closed lifted uncertainty}
\begin{align}
\minimize \quad & x_1 \\
\st \quad & \bm{\xi}^{\top} \left( \bm{A_n} \bm{x} + \bm{D_n} \bm{y} \right) + \bm{\tilde{a}_n}^{\top} \bm{\tilde{x}}(\bm{\xi'}) + \bm{\tilde{d}_n}^{\top} \bm{\tilde{y}}(\bm{\xi'}) \leq \bm{\xi}^{\top} \bm{b_n} \quad \forall \, n \in \mathcal{N}, \, \left( \bm{\xi}, \, \bm{\xi'} \right) \in \dot{\widetilde{\Xi}} \left( \bm{y} \right) \label{eqn:GeneralConstraintLiftedUncertainty} \\
& \bm{x} \in \mathbb{R}^{P}, \; \bm{y} \in \{0,1\}^{Q} \\
& \bm{\tilde{x}}(\bm{\xi'}) \in \mathbb{R}^{\widetilde{P}}, \; \bm{\tilde{y}}(\bm{\xi'}) \in \{0,1\}^{\widetilde{Q}} \quad \forall \, \bm{\xi'} \in \dot{\widetilde{\Xi}}_{\bm{\xi}'}(\bm{y}), \label{eqn:VariableDomainLiftedUncertainty}
\end{align}
\end{subequations}
which has the same optimal value as \eqref{eqn:TwoStage}, and there is a one-to-one mapping between feasible and optimal solutions of problems \eqref{eqn:TwoStage} and \eqref{eqn:TwoStage with closed lifted uncertainty}, as shown in \citet{Bertsimas2018}.

However, because the uncertain parameters considered are correlated and decision-dependent, it is difficult to formulate the exact closed-form expression for $\dot{\widetilde{\Xi}}(\bm{y})$. Hence, we introduce a tractable outer approximation $\widetilde{\Xi}(\bm{y})$, which stems from the convex hull of the marginal support of $\bm{\tilde{\xi}}_i$ for all $i = 1,\dots,K$. By replacing $\dot{\widetilde{\Xi}}(\bm{y})$ with $\widetilde{\Xi}(\bm{y})$, the two-stage problem is subsequently transformed into \eqref{eqn:TwoStage with polyhedral closed lifted uncertainty}, whose solution is rather a conservative approximation of the one derived by \eqref{eqn:TwoStage with closed lifted uncertainty}:
\begin{subequations}
\label{eqn:TwoStage with polyhedral closed lifted uncertainty}
\begin{align}
\minimize \quad & x_1 \\
\st \quad & \bm{\xi}^{\top} \left( \bm{A_n} \bm{x} + \bm{D_n} \bm{y} \right) + \bm{\tilde{a}_n}^{\top} \bm{\tilde{x}}(\bm{\xi'}) + \bm{\tilde{d}_n}^{\top} \bm{\tilde{y}}(\bm{\xi'}) \leq \bm{\xi}^{\top} \bm{b_n} \quad \forall \, n \in \mathcal{N}, \, \left( \bm{\xi}, \, \bm{\xi'} \right) \in \widetilde{\Xi} \left( \bm{y} \right) \label{eqn:PolyhedralGeneralConstraintLiftedUncertainty} \\
& \bm{x} \in \mathbb{R}^{P}, \; \bm{y} \in \{0,1\}^{Q} \\
& \bm{\tilde{x}}(\bm{\xi'}) \in \mathbb{R}^{\widetilde{P}}, \; \bm{\tilde{y}}(\bm{\xi'}) \in \{0,1\}^{\widetilde{Q}} \quad \forall \, \bm{\xi'} \in \widetilde{\Xi}_{\bm{\xi}'}(\bm{y}). \label{eqn:PolyhedralVariableDomainLiftedUncertainty}
\end{align}
\end{subequations}

To obtain $\widetilde{\Xi}(\bm{y})$, first notice that both $\bar{L}$ and $\widehat{L}$ can be considered piecewise mapping functions with each piece being a subinterval $\Xi_i^j$ formed by two consecutive breakpoints:
\begin{equation*} \label{eqn: subinterval}
    \Xi_i^j = \left\{ \xi_i \in \mathbb{R}:\;  p_i^{j-1} \leq \xi_i \leq p_i^j \right\}.
\end{equation*}
For ease of exposition, we set $p_i^0 = \xi_i^{\min}$ and $p_i^{r_i} = \xi_i^{\max}$. Let $\mathcal{V}^j_i = \{p_i^{j-1}, p_i^j\}$, then for every $i = 1,\ldots,K$, we can define the following sets of vertices in the lifted space:
\begin{equation*}
%\begin{aligned}
\begin{split}
     &\widetilde{\mathcal{V}}_i  =  \bigcup_{j=1}^{r_i} \widetilde{\mathcal{V}}_{i}^j \\
     &\widetilde{\mathcal{V}}_i^j  = \left \{
     \bm{\tilde{v}}_i =  \left (v_i, \, \bm{\bar{v}}_i, \, \bm{\hat{v}}_i \right): 
     \begin{array}{l}
          v_i \in \mathcal{V}_i^j \\
         \bm{\bar{v}}_i= \bar{L}_i (\bm{e}_i v_i) \\ 
         \bm{\hat{v}}_i = \lim\limits_{\xi_i \rightarrow v_i, \, \xi_i \in \Xi_i^j} \widehat{L}_i (\bm{e}_i \xi_i) 
     \end{array}
    \right \} \quad \forall \, j = 1, \dots, r_i,
\end{split}{}
%\end{aligned}
\end{equation*}
which allow us to formulate the following convex hull representation for the closure of the marginal support of $\bm{\tilde{\xi}}_i$:
\begin{equation} \label{eqn: LiftedSubset}
   \widetilde{\Xi}_i = \left \{
    \bm{\tilde{\xi}}_i = \left (\xi_i, \bm{\bar{\xi}}_i, \bm{\hat{\xi}}_i \right):
        \begin{array}{l}
            %  \exists \lambda_i^j(\bm{\tilde{v}}_i) \in \mathbb{R}_+, \; \forall \bm{\tilde{v}}_i \in \widetilde{V}_i^j, \; j=1, \ldots, r_i \; \text{such that} \\
             \sum\limits_{j=1}^{r_i} \sum\limits_{\bm{\tilde{v}}_i \in \widetilde{\mathcal{V}}_{i}^j} \lambda_i^j(\bm{\tilde{v}}_i) = 1 \\
             \xi_i = \sum\limits_{j=1}^{r_i} \sum\limits_{\bm{\tilde{v}}_i \in \widetilde{\mathcal{V}}_{i}^j} \lambda_i^j(\bm{\tilde{v}}_i) v_i \\
             \bm{\bar{\xi}}_i = \sum\limits_{j=1}^{r_i} \sum\limits_{\bm{\tilde{v}}_i \in \widetilde{\mathcal{V}}_{i}^j} \lambda_i^j(\bm{\tilde{v}}_i) \bm{\bar{v}}_i\\
             \bm{\hat{\xi}}_i = \sum\limits_{j=1}^{r_i} \sum\limits_{\bm{\tilde{v}}_i \in \widetilde{\mathcal{V}}_{i}^j} \lambda_i^j(\bm{\tilde{v}}_i) \bm{\hat{v}}_i  \\
             \lambda_i^j(\bm{\tilde{v}}_i) \in \mathbb{R}_+ \quad \forall \, j=1, \ldots, r_i, \; \bm{\tilde{v}}_i \in \widetilde{\mathcal{V}}_i^j
            \end{array}
            \right \},
\end{equation}
where $\lambda_i^j(\bm{\tilde{v}}_i)$ denotes the coefficient associated with a particular vertex $\bm{\tilde{v}}_i \in \widetilde{\mathcal{V}}_{i}^j$. If the original uncertain parameters are independent and exogenous, $\dot{\widetilde{\Xi}}$ can be exactly represented as the Cartesian product of all $\widetilde{\Xi}_i$. In the general and endogenous case, we have
\begin{equation} \label{eqn:AbstractSet}
\dot{\widetilde{\Xi}}(\bm{y}) \subseteq \widetilde{\Xi}(\bm{y}) := \left \{ \prod_{i=1}^K \widetilde{\Xi}_i \right \} \bigcap \left \lbrace
     \left (\bm{\xi}, \, \bm{\xi}' \right) \in \mathbb{R}^{K + K'}:
  \bm{W} \bm{\xi} \leq \bm{U} \bm{y}, \; \xi_1 = 1 \right \rbrace.
\end{equation}
Since $\widetilde{\Xi}(\bm{y})$ is generally a superset of $\dot{\widetilde{\Xi}}(\bm{y})$, it may result in a more conservative solution. However, it is worth mentioning that the uncertainty set in the space of the original uncertain parameters, i.e. the projection of $\widetilde{\Xi}(\bm{y})$ onto the $\bm{\xi}$-space, remains unchanged. This outer approximation only applies to the new lifted uncertain parameters $\bm{\xi}'$, which are used for the construction of the decision rules, as we will show in Section \ref{sec:DecisionRules}. Therefore, the use of this outer approximation should not be interpreted as expanding the uncertainty set to robustify against, but rather as further restricting the set of possible decision rules. As such, this increased conservatism could be compensated by adjusting the locations or increasing the number of breakpoints.

\begin{example}
In this small example, we illustrate the relationship between $\dot{\widetilde{\Xi}}(y)$ and $\widetilde{\Xi}(y)$. Consider the following decision-dependent uncertainty set:
\begin{equation*}
    \Xi(y) = \left \lbrace \bm{\xi} \in \mathbb{R}^2: \, 0 \leq \xi_2 \leq 5 \xi_1 - y, \; \xi_1  = 1 \right \rbrace.
\end{equation*}
We place one breakpoint, $p_2^1$, at the center of the marginal support of $\xi_2$, and analyze the resulting $\dot{\widetilde{\Xi}}(y)$ and $\widetilde{\Xi}(y)$. For ease of visualization, we only show the projections of these two sets onto the two-dimensional $(\xi_2,\hat{\xi}_2^1)$-space for different values of $y$. If $y=0$, $\dot{\widetilde{\Xi}}(y)$ and $\widetilde{\Xi}(y)$ are in fact the same, which is depicted by the gray-shaded area in Figure \ref{fig: OuterApprox-Original}. However, if $y = 1$, the projection of $\dot{\widetilde{\Xi}}(y)$ is the red-shaded area in Figure \ref{fig: OuterApprox-True}, while that of $\widetilde{\Xi}(y)$ is the green-shaded area in Figure \ref{fig: OuterApprox-Outer}. Clearly, $\dot{\widetilde{\Xi}}(y)$ is a subset of $\widetilde{\Xi}(y)$.

\begin{figure}[ht]\centering
\subfloat[$\dot{\widetilde{\Xi}}(\bm{y}) = \widetilde{\Xi}(\bm{y})$ for $y = 0$]{
  \label{fig: OuterApprox-Original}
	\fbox{\includegraphics[height=1.2in]{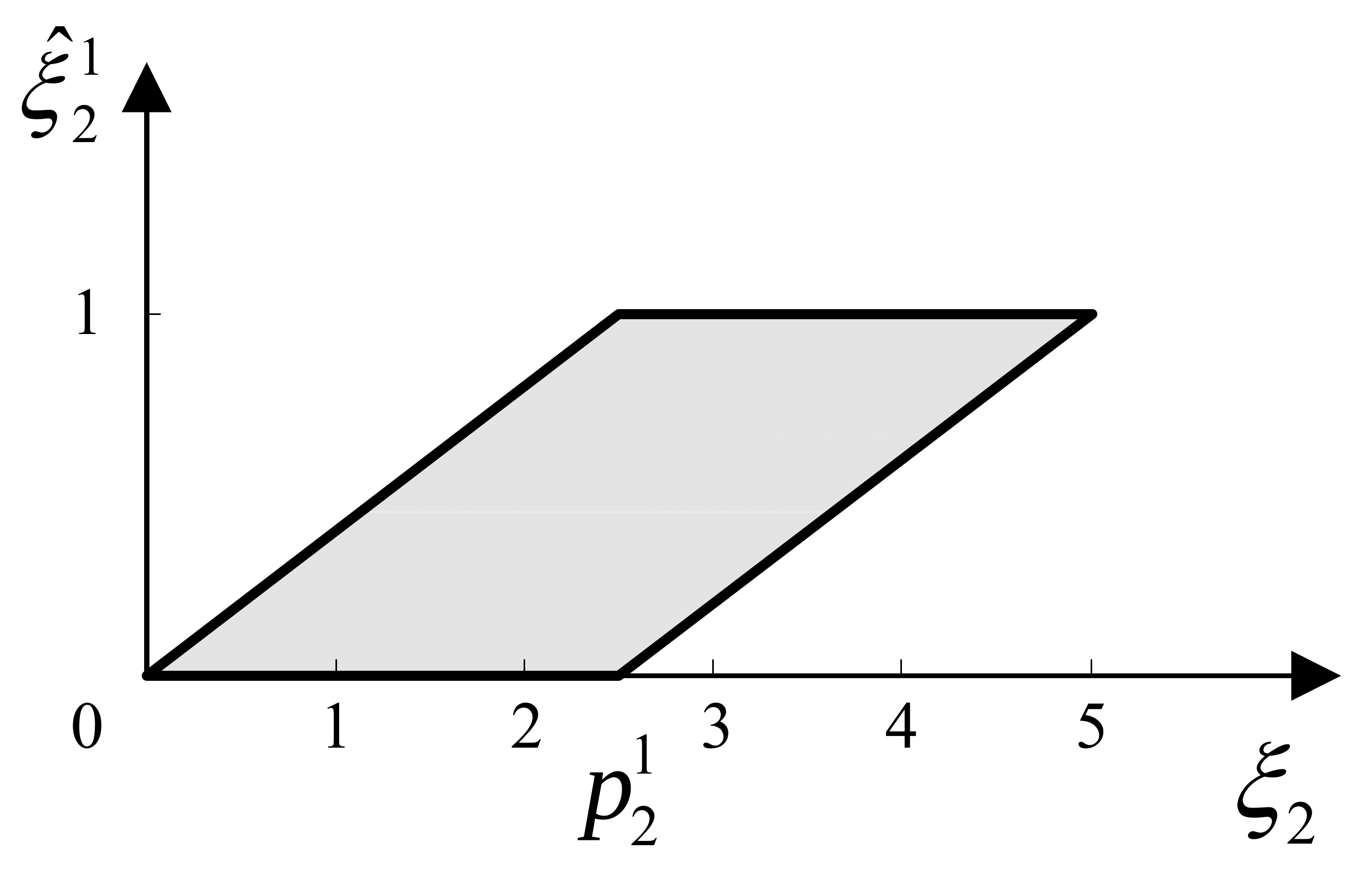}}}
\subfloat[$\dot{\widetilde{\Xi}}(\bm{y})$ for $y = 1$]{
  \label{fig: OuterApprox-True}
	\fbox{\includegraphics[height=1.2in]{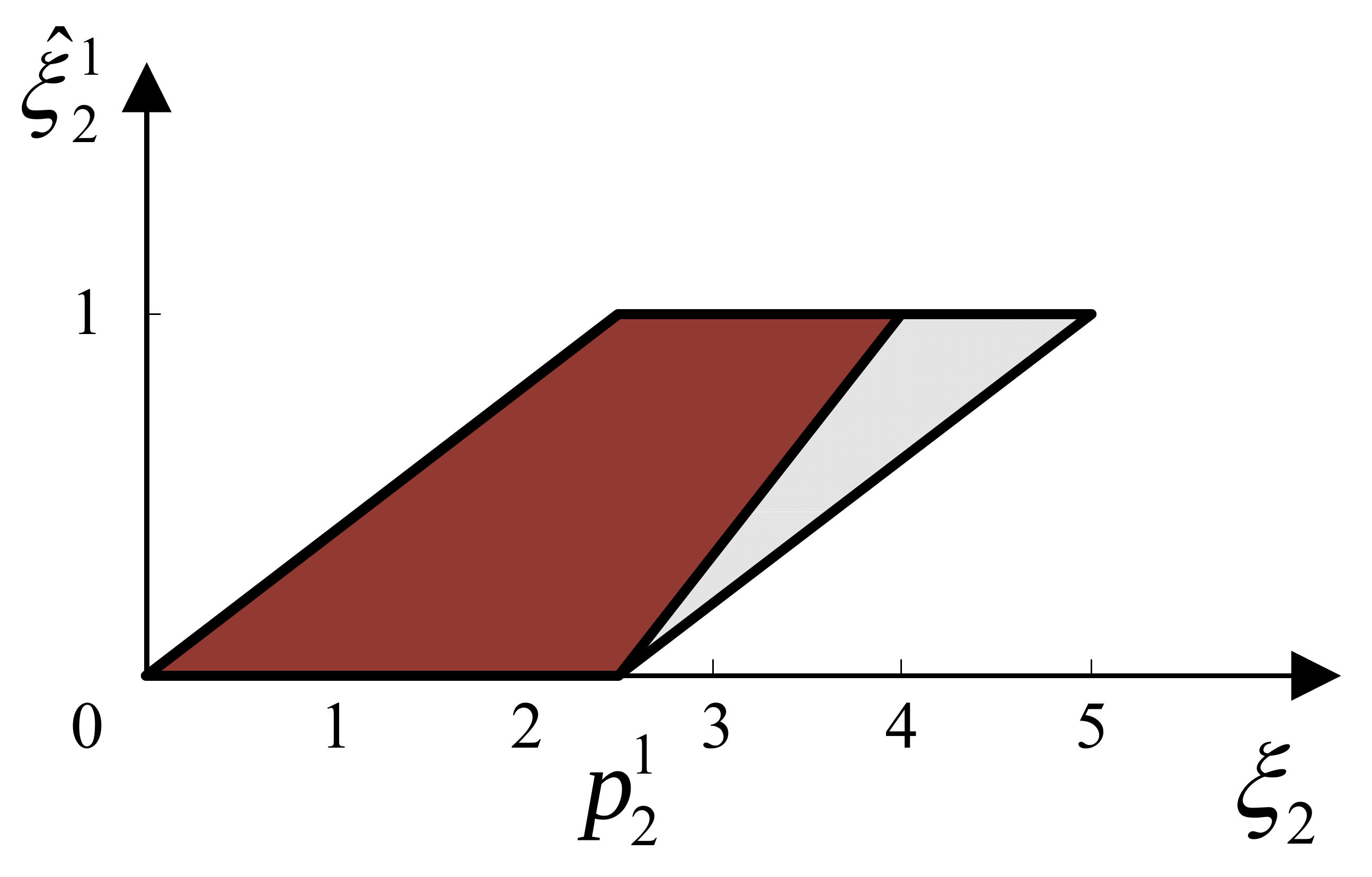}}}
\subfloat[$\widetilde{\Xi}(\bm{y})$ for $y = 1$]{
  \label{fig: OuterApprox-Outer}
	\fbox{\includegraphics[height=1.2in]{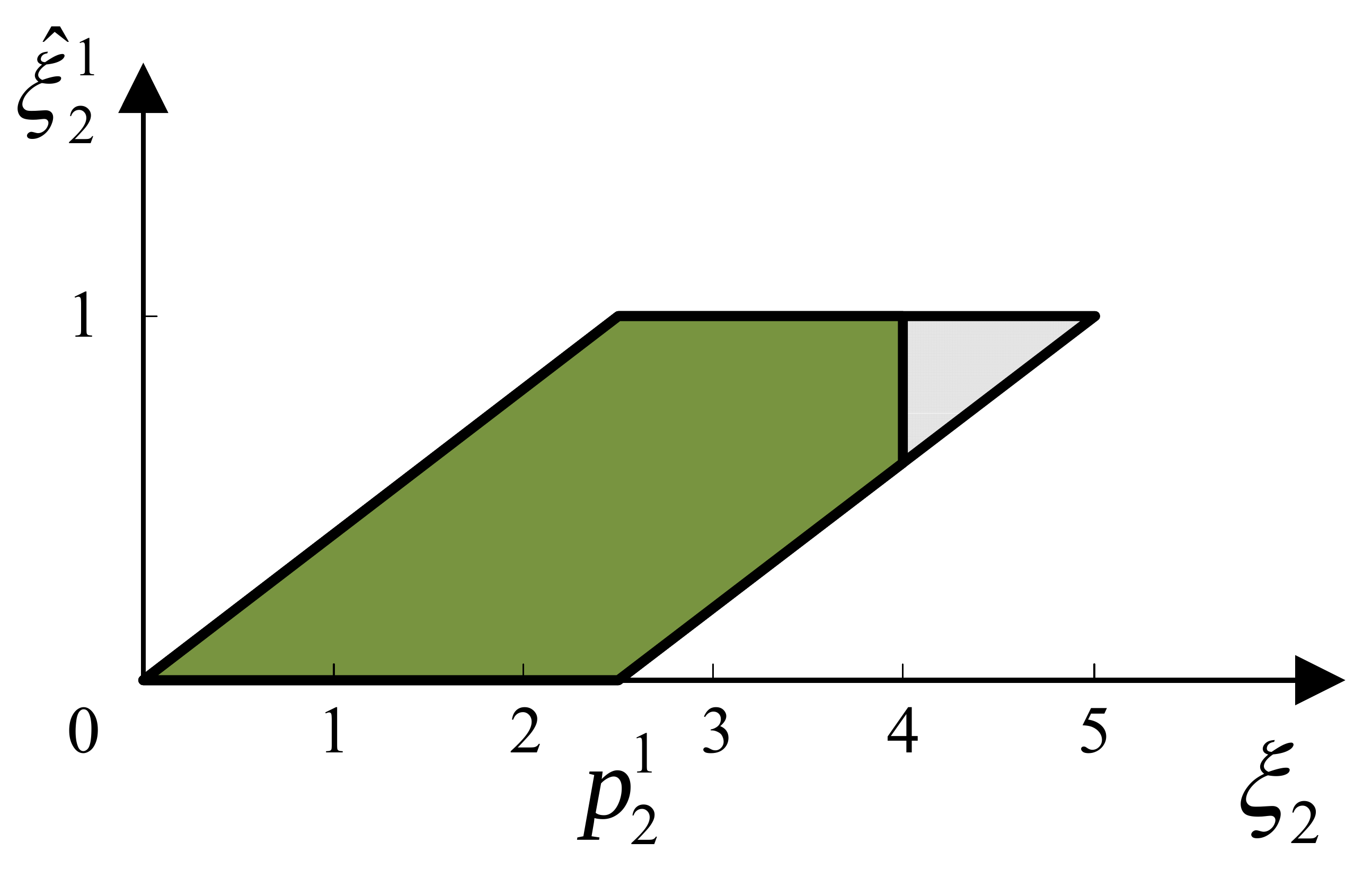}}}
\caption{Projections of $\dot{\widetilde{\Xi}}(y)$ and $\widetilde{\Xi}(y)$ onto the $(\xi_2,\hat{\xi}_2^1)$-space for $y=0$ and $y=1$.}
\label{fig:OuterApprox}
\end{figure}
\end{example}

\subsection{Decision Rule Approximation and Reformulation}
\label{sec:DecisionRules}

Following the decision rule approach, we solve problem \eqref{eqn:TwoStage with polyhedral closed lifted uncertainty} approximately by restricting the recourse decisions to adjust according to some decision rules, which are functions of the uncertain parameters. Here, we apply decision rules parameterized as follows:
\begin{subequations}
\label{eqn:ProposedLinearRule}
\begin{align}
    & \bm{\tilde{x}}=\sum_{i=1}^K \left( \bm{\overline{X}}_i \bm{\bar{\xi}}_i + \bm{\widehat{X}}_i \bm{\hat{\xi}}_i \right)
    \label{eqn: continuous decision rule}\\
    & \bm{\tilde{y}} = \sum_{i=1}^K \bm{\widehat{Y}}_i \bm{\hat{\xi}}_i
    \label{eqn: binary decision rule}
\end{align}
\end{subequations}
where $\bm{\overline{X}}_i \in \mathbb{R}^{\widetilde{P} \times r_i}$, $\bm{\widehat{X}}_i \in \mathbb{R}^{\widetilde{P} \times g_i}$, and $\bm{\widehat{Y}}_i \in \{-1, 0, 1\}^{\widetilde{Q} \times g_i}$ for $i=1, \ldots, K$. We emphasize that the structure of \eqref{eqn: continuous decision rule} allows \textit{discontinuous} piecewise linear decision rules; this is in contrast to most existing works in the literature, which use affine or continuous piecewise linear decisions rules for continuous recourse variables. Also note that, as shown in \citet{Bertsimas2018}, the given domain for $\bm{\widehat{Y}}_i$ is sufficient to allow binary $\bm{\tilde{y}}$ using the decision rule in \eqref{eqn: binary decision rule}. Before proceeding to the reformulation, we illustrate with the following simple example that, unlike in continuous optimization, it is often crucial for the continuous recourse variables to follow discontinuous decision rules in mixed-integer optimization.

\begin{example}
Consider the following MILP:
\begin{equation}
\label{eqn:DecisionRuleExample}
\begin{aligned}
    \minimize \quad & 2\tilde{x}_1 +  \tilde{x}_2 \\
    \st \quad & \tilde{y}_1 \leq \tilde{x}_1 \leq 3\tilde{y}_1 \\
    & 3\tilde{y}_2 \leq \tilde{x}_2 \leq 5\tilde{y}_2 \\
    & \tilde{x}_1 + \tilde{x}_2 = \xi \\
    & \tilde{y}_1 + \tilde{y}_2 \leq 1 \\
    & \tilde{x}_1, \, \tilde{x}_2 \in \mathbb{R} \\
    & \tilde{y}_1, \, \tilde{y}_2 \in \left \{0,1 \right \},
    \end{aligned}
\end{equation}
where the continuous variables $\tilde{x}_1$ and $\tilde{x}_2$ take nonzero values if and only if the respective binary variables $\tilde{y}_1$ and $\tilde{y}_2$ are equal to 1. For $\xi \in \Xi = \left \{ \xi \in \mathbb{R}: \; 1 \leq \xi \leq 5\right \}$, the optimal solution to \eqref{eqn:DecisionRuleExample} as a function of $\xi$ is as follows:
\begin{equation}
\label{eqn:DecisionRuleExampleSolution}
\left \{
    \begin{array}{ll}
     \tilde{y}_1 = 1, \; \tilde{x}_1 = \xi , \; \tilde{y}_2 = \tilde{x}_2 = 0  \quad  &  \text{if} \; \, 1 \leq \xi < 3  \\
     \tilde{y}_1 = \tilde{x}_1 = 0, \; \tilde{y}_2 = 1, \; \tilde{x}_2 = \xi \quad   &   \text{if} \; \, 3 \leq \xi \leq 5
    \end{array}{}
    \right .
\end{equation}
which is shown in Figure \ref{fig:EX2}. One can see that $\tilde{x}_1$ and $\tilde{x}_2$ follow discontinuous piecewise linear functions. In fact, \eqref{eqn:DecisionRuleExampleSolution} represents the only feasible solution when $\xi \neq 3$ and the only optimal solution when $\xi = 3$. As a result, restricting the continuous variables to be continuous piecewise linear functions of $\xi$ would render the problem infeasible.

\begin{figure}[ht]\centering
    \subfloat[$x_1$ and $y_1$]{
        \label{fig:EX2-1}
        \fbox{\includegraphics[width = 2.5in]{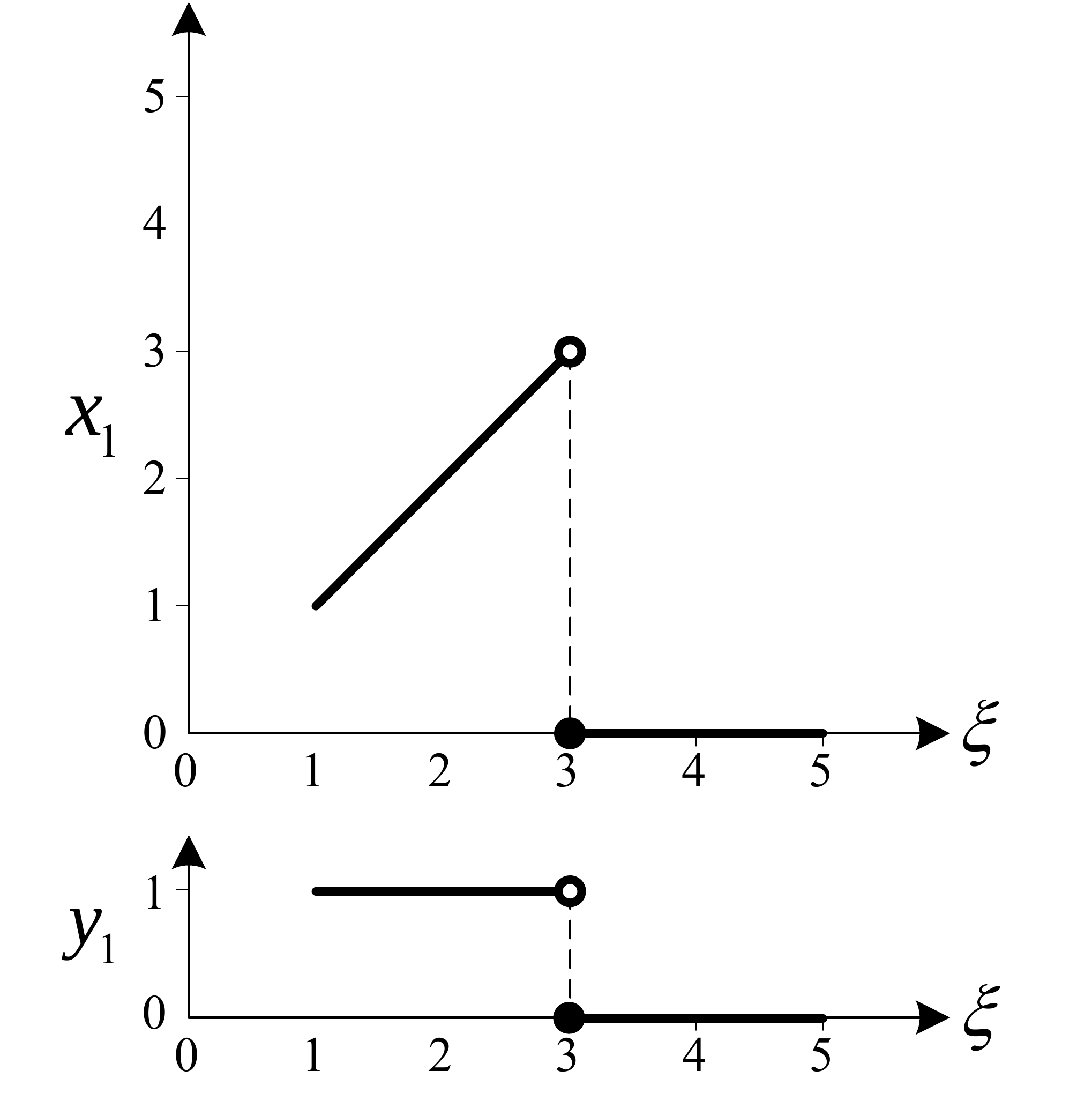}}}
    \subfloat[$x_2$ and $y_2$]{
        \label{fig:EX2-2}
        \fbox{\includegraphics[width = 2.5in]{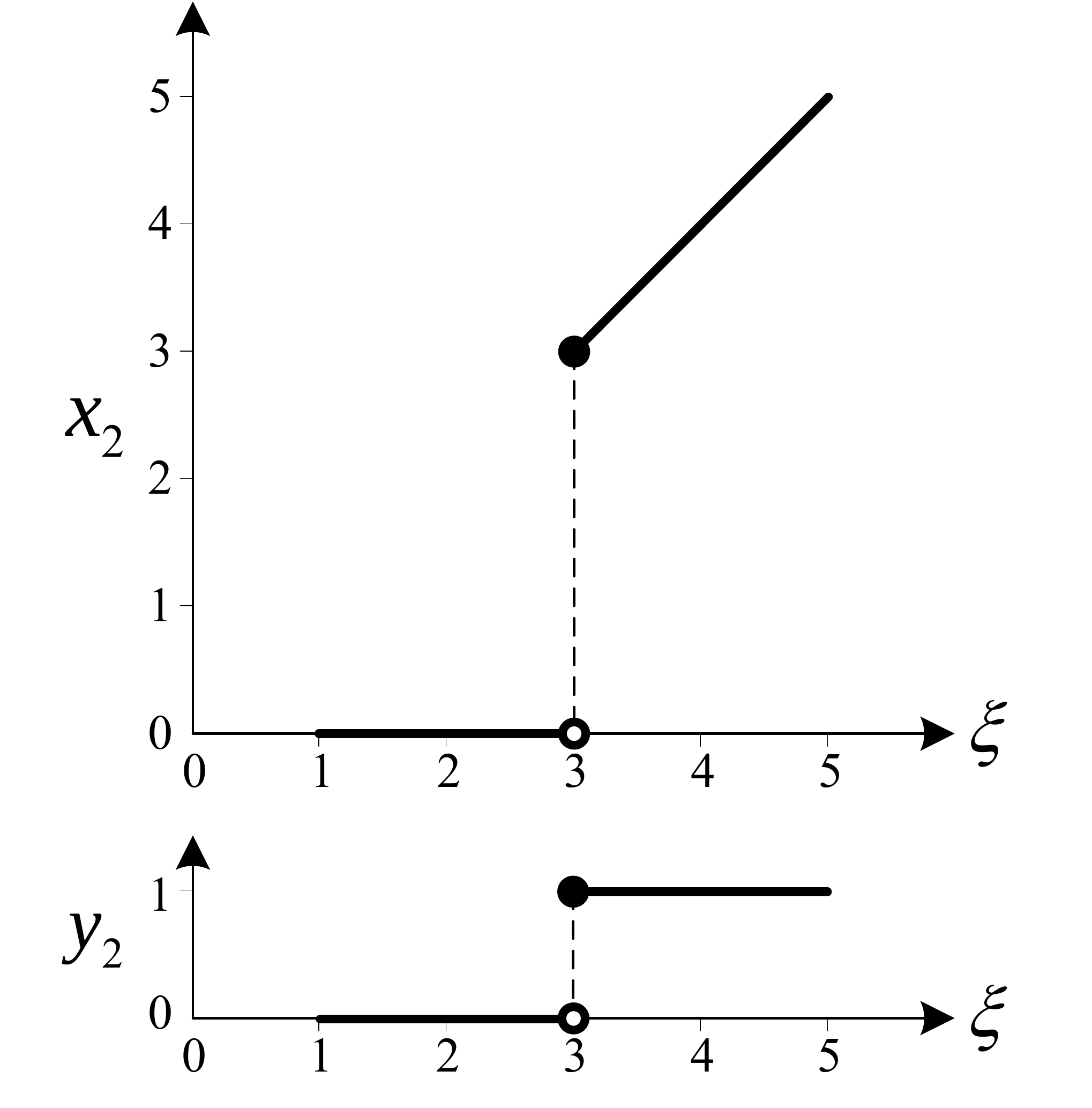}}}
\caption{Continuous/binary variables as discontinuous piecewise linear/constant functions of $\xi$.}
\label{fig:EX2}
\end{figure} 
\end{example}

By substituting the decisions rules \eqref{eqn:ProposedLinearRule} into \eqref{eqn:TwoStage with polyhedral closed lifted uncertainty}, constraints \eqref{eqn:PolyhedralGeneralConstraintLiftedUncertainty} become
\begin{equation}
\label{eqn:PreliminaryReformulation}
    \bm{\xi}^{\top} \bm{f}_n + \sum_{i=1}^K \left[
    \bm{\tilde{a}}_n^{\top} \bm{\overline{X}}_i \bm{\bar{\xi}}_i + 
    \left (
    \bm{\tilde{a}}_n^{\top} \bm{\widehat{X}}_i + \bm{\tilde{d}}_n^{\top}\widehat{Y}_i
    \right ) \bm{\hat{\xi}}_i \right]
    \leq 0 \quad \forall \, n \in \mathcal{N}, \; \left (\bm{\xi}, \, \bm{\bar{\xi}}, \, \bm{\hat{\xi}} \right) \in \widetilde{\Xi}(\bm{y}) 
\end{equation}
with $\bm{f}_n = \bm{A}_n \bm{x} + \bm{D}_n \bm{y} - \bm{b}_n$. Following standard robust optimization techniques, we first apply the worst-case reformulation:
\begin{equation}
    \label{eqn:Worst-caseReformulation}
    \max_{\left (\bm{\xi}, \, \bm{\bar{\xi}}, \, \bm{\hat{\xi}} \right) \in \widetilde{\Xi}(\bm{y}) } \left\lbrace \bm{\xi}^{\top} \bm{f}_n + \sum_{i=1}^K \left[
    \bm{\tilde{a}}_n^{\top} \bm{\overline{X}}_i \bm{\bar{\xi}}_i + 
    \left (
    \bm{\tilde{a}}_n^{\top} \bm{\widehat{X}}_i + \bm{\tilde{d}}_n^{\top}\widehat{Y}_i
    \right ) \bm{\hat{\xi}}_i \right] \right\rbrace
    \leq 0 \quad \forall \, n \in \mathcal{N},
\end{equation}
which then leads to the following reformulation due to strong duality of the left-hand-side maximization problems:
\begin{equation} \label{eqn: reformulation of 1st constraint}
\begin{aligned}
    & \bm{e}^{\top} \bm{\rho}_n + \bm{\mu}_n^{\top} \bm{U} \bm{y} \leq 0 \quad \forall \, n \in \mathcal{N} \\
    &  \rho_{ni} + \bm{\mu}_n^{\top} \bm{w}_i v_i \geq f_{ni} v_i + \bm{\tilde{a}}_n^{\top} \bm{\overline{X}}_i \bm{\bar{v}}_i + \left( \bm{\tilde{a}}_n^{\top} \bm{\widehat{X}}_i + \bm{\tilde{d}}_n^{\top} \bm{\widehat{Y}}_i \right) \bm{\hat{v}}_i \\
    & \quad \quad \quad \quad \quad \quad \quad \quad \quad \quad \forall \, n \in \mathcal{N}, \; i=1, \ldots, K, \;  \left(v_i, \, \bm{\bar{v}}_i, \, \bm{\hat{v}}_i \right) \in \widetilde{\mathcal{V}}_i \\
   &  \bm{\rho}_n \in \mathbb{R}^K, \; \bm{\mu}_n \in \mathbb{R}_+^M \quad \forall \, n \in \mathcal{N},
\end{aligned}
\end{equation}
where $\bm{\rho}_n$ and $\bm{\mu}_n$ are the dual variables of the maximization problem,  $f_{ni}$ and $\rho_{ni}$ denote the $i$th element of $\bm{f}_n$ and $\bm{\rho}_n$ respectively, and $\bm{w}_i$ is the $i$th column vector of matrix $\bm{W}$.

In addition, since the decision rule in \eqref{eqn: binary decision rule} is guaranteed to yield integer $\bm{\tilde{y}}$, the integratility constraints $\bm{\tilde{y}}(\bm{\xi'}) \in \left \{ 0,1 \right \}^{\widetilde{Q}}$ in \eqref{eqn:PolyhedralVariableDomainLiftedUncertainty} can be relaxed to
\begin{equation*}
    \bm{0} \leq \bm{\tilde{y}}(\bm{\xi'}) \leq \bm{e} \quad \forall \, \bm{\xi'} \in \widetilde{\Xi}_{\bm{\xi}'}(\bm{y}),
\end{equation*}
which, using similar arguments as above, can be reformulated into the following set of constraints:
\begin{equation}
\label{eqn: reformulation of 2nd constraint}
\begin{aligned}
        &  \bm{\underline{\Omega}} \bm{e} - \bm{\underline{\Psi}} \bm{U} \bm{y} \geq \bm{0} \\
         & \bm{\widehat{Y}}_i \bm{\hat{v}}_i - \bm{\underline{\omega}}_i + \bm{\underline{\Psi}} \bm{w}_i v_i \geq \bm{0} \quad \forall \, i=1, \ldots, K,  \; \left(v_i, \, \bm{\bar{v}}_i, \, \bm{\hat{v}}_i \right) \in \widetilde{\mathcal{V}}_i \\
        & \bm{\overline{\Omega}} \bm{e} + \bm{\overline{\Psi}} \bm{U} \bm{y} \leq \bm{e} \\
        & \bm{\widehat{Y}}_i \bm{\hat{v}}_i - \bm{\overline{\omega}}_i - \bm{\overline{\Psi}} \bm{w}_i v_i \leq \bm{0} \quad \forall \, i=1, \ldots, K, \;  \left(v_i, \, \bm{\bar{v}}_i, \, \bm{\hat{v}}_i \right) \in \widetilde{\mathcal{V}}_i \\
        & \bm{\underline{\Omega}}, \, \bm{\overline{\Omega}} \in \mathbb{R}^{\widetilde{Q} \times K}, \quad \bm{\underline{\Psi}}, \, \bm{\overline{\Psi}} \in \mathbb{R}_+^{\widetilde{Q} \times M},
        \end{aligned}
\end{equation}
where $\bm{\underline{\Omega}}$, $\bm{\overline{\Omega}}$, $\bm{\underline{\Psi}}$ and $\bm{\overline{\Psi}}$ are the matrices of dual variables, and $\bm{\underline{\omega}}_i$ and $\bm{\overline{\omega}}_i$ are the $i$th columns of $\bm{\underline{\Omega}}$ and $\bm{\overline{\Omega}}$, respectively.

Finally, we arrive at the following formulation:
\begin{equation}
\label{eqn: two-stage reformulation}
\begin{aligned}
    \minimize \quad & x_1 \\
    \st \quad & \bm{f}_n = \bm{A}_n \bm{x} + \bm{D}_n \bm{y}- \bm{b}_n \quad \forall \, n \in \mathcal{N} \\
    & \eqref{eqn: reformulation of 1st constraint}, \; \eqref{eqn: reformulation of 2nd constraint} \\
    & \bm{x} \in \mathbb{R}^{P}, \; \bm{y} \in \{0,1\}^{Q} \\
    & \bm{\overline{X}}_i \in \mathbb{R}^{\widetilde{P} \times r_i}, \; \bm{\widehat{X}}_i \in \mathbb{R}^{\widetilde{P} \times g_i}, \; \bm{\widehat{Y}}_i \in \{-1, 0, 1\}^{\widetilde{Q} \times g_i} \quad \forall \, i=1, \ldots, K.
\end{aligned}
\end{equation}
Note that \eqref{eqn: reformulation of 1st constraint} and \eqref{eqn: reformulation of 2nd constraint} contain bilinear terms; however, each bilinear term consists of a continuous variable multiplied by a binary variable and can hence be exactly linearized \citep{Glover1975}. As a result, problem \eqref{eqn: two-stage reformulation} can be solved as an MILP.

\section{The Multistage Case}
\label{sec:Multistage}

In this section, the proposed methodology is extended to multistage robust MILPs with endogenous uncertainty of the following form:
\begin{subequations}
\label{model:GeneralMultistageDeterministicModel}
\begin{align}
    \minimize
    % _{\bm{x}_1, \bm{y}_1, \bm{x}_2(\cdot),\ldots, \bm{x}_T(\cdot), \bm{y}_2(\cdot),\ldots, \bm{y}_T(\cdot)} 
    \quad & \bm{e}_1^{\top} \bm{x}_1 \\
    \st \quad & \bm{A}^1 \bm{x}_1 + \bm{D}^1 \bm{y}_1 \leq \bm{b}^1 \\
    & \bm{A}^t(\bm{\xi}^{[t]}) \bm{x}_1 + \bm{D}^t(\bm{\xi}^{[t]}) \bm{y}_1 + \sum\limits_{t'=2}^t \left[ \bm{\widetilde{A}}_{t'}^t \bm{x}_{t'}(\bm{\xi}^{[t']}) + \bm{\widetilde{D}}_{t'}^t \bm{y}_{t'}(\bm{\xi}^{[t']}) \right] \leq \bm{b}^t (\bm{\xi}^{[t]}) \notag \\
    & \quad\quad\quad\quad\quad\quad\quad\quad\quad\quad\quad\quad\quad\quad\quad \forall \, t = 2, \ldots, T, \; \bm{\xi}^{[t]} \in \Xi^{[t]}(\bm{y}^{[t-1]}) \label{eqn:MultiStageDeterministicConstraint}\\
    & \bm{x}_1 \in \mathbb{R}^{P_1}, \; \bm{y}_1 \in \{ 0, 1 \}^{Q_1} \\
    & \bm{x}_t(\bm{\xi}^{[t]}) \in \mathbb{R}^{P_t}, \; \bm{y}_t(\bm{\xi}^{[t]}) \in \{ 0, 1 \}^{Q_t} \quad \forall \, t = 2, \ldots, T, \; \bm{\xi}^{[t]} \in \Xi^{[t]}( \bm{y}^{[t-1]}) \label{eqn:MultiStageDomainConstraint}
\end{align}
\end{subequations}
where $\bm{x}_t$ and $\bm{y}_t$ are the continuous and binary variables, respectively, in stage $t$. The uncertain parameters whose true values are observed in stage $t$ are denoted by $\bm{\xi}_t = (\xi_{t1}, \ldots, \xi_{tK_t}) \in \mathbb{R}^{K_t}$. Again for notational convenience, $K_1 = 1$ and the only element of $\bm{\xi}_1$, $\xi_{11}$, is assumed to be 1. We further define all uncertain parameters observed up to stage $t$ as $\bm{\xi}^{[t]} = \left( \bm{\xi}_1, \bm{\xi}_2, \ldots, \bm{\xi}_t \right)$. Moreover, $\bm{A}^t(\bm{\xi}^{[t]}) \in \mathbb{R}^{N_t \times P_1}$, $\bm{D}^t(\bm{\xi}^{[t]}) \in \mathbb{R}^{N_t \times Q_1}$, and $\bm{b}^t (\bm{\xi}^{[t]}) \in \mathbb{R}^{N_t}$ are coefficient matrices or vectors that depend linearly on $\bm{\xi}^{[t]}$, i.e.
\begin{equation}
\label{eqn:LMIinFirstStage}
     \bm{A}^t(\bm{\xi}^{[t]}) = \sum\limits_{t'=1}^t \sum\limits_{i=1}^{K_{t'}} \bm{A}_{t'i}^t \xi_{t'i}, \quad
     \bm{D}^t(\bm{\xi}^{[t]}) = \sum\limits_{t'=1}^t \sum\limits_{i=1}^{K_{t'}} \bm{D}_{t'i}^t \xi_{t'i}, \quad 
     \bm{b}^t (\bm{\xi}^{[t]}) = \sum\limits_{t'=1}^t \sum\limits_{i=1}^{K_{t'}} \bm{b}_{t'i}^t \xi_{t'i},
\end{equation}
where $\bm{A}_{t'i}^t \in \mathbb{R}^{N_t \times P_1}$, $\bm{D}_{t'i}^t \in \mathbb{R}^{N_t \times Q_1}$, and $\bm{b}_{t'i}^t \in \mathbb{R}^{N_t}$.

Formulation \eqref{model:GeneralMultistageDeterministicModel} indicates that the uncertainty set $\Xi^{[t]}$ in stage $t$ depends on binary decisions made in previous stages, $\bm{y}^{[t-1]} = \left(\bm{y}_1, \dots, \bm{y}_{t-1}\right)$, where $\bm{y}_{\hat{t}}(\bm{\xi}^{[\hat{t}]})$ are recourse variables for $\hat{t} \geq 2$. This means that, in contrast to the two-stage case, the changing of the uncertainty set is now also a recourse decision. Here, we consider decision-dependent uncertainty sets of the following form:
\begin{equation}
\label{eqn:MultiStageOriginalUncertaintySet}
    \Xi^{[t]}(\bm{y}^{[t-1]}) = \left \{
    \bm{\xi}^{[t]} \in \mathbb{R}^{K^{[t]}}:
    \bm{W}^{[t]} \bm{\xi}^{[t]} \leq  \sum\limits_{\hat{t}=1}^{t-1} \bm{U}_{\hat{t}}^t \bm{y}_{\hat{t}}(\bm{\xi}^{[\hat{t}]}), \; \xi_{11} = 1
    \right \},
\end{equation}
where $K^{[t]} = \sum_{t'=1}^t K_{t'}$, $\bm{W}^{[t]} \in \mathbb{R}^{M_t \times K^{[t]}}$ and $\bm{U}_{\hat{t}}^t \in \mathbb{R}^{M_t \times Q_{\hat{t}}}$. We assume that $\Xi^{[t]}(\bm{y}^{[t-1]})$ is a compact polyhedron for any feasible $\bm{y}^{[t-1]}$. Note that $\Xi^{[t]}(\bm{y}^{[t-1]})$ reduces to the uncertainty set in the two-stage problem shown in \eqref{eqn:UncertaintySet} when $t=2$.

\subsection{Lifted Uncertainty Set}
% With $(r_{t'i}-1)$ breakpoints placed inside the marginal support of each uncertain parameter $\xi_{t'i}$, the two previously defined lifting operators can also be adapted to the multi-stage case, which are $\bar{L}: \mathbb{R}^{K^{[t]}} \mapsto \mathbb{R}^{\bar{K}^{[t]}}$ and $\widehat{L}: \mathbb{R}^{{K}^{[t]}} \mapsto \mathbb{R}^{\widehat{K}^{[t]}}$ with $\bar{K}^{[t]} = \sum_{t'=1}^t \sum_{i=1}^{K_{t'}} r_{t'i}$ and $\widehat{K}^{[t]} = \sum_{t'=1}^t \sum_{i=1}^{K_{t'}} g_{t'i}$, where $g_{t'i} = \max \left \{ 1, r_{t'i}-1 \right \}$. Then, the vectors of lifted uncertainty from $\bm{\xi}^{[t]}$ by $\bar{L}$ and $\widehat{L}$ are denoted by $\bm{\bar{\xi}}^{[t]} = \left (\bm{\bar{\xi}}_{11}, \bm{\bar{\xi}}_{21}, \ldots, \bm{\bar{\xi}}_{t K_t}  \right)$ and $\bm{\hat{\xi}}^{[t]} = \left (\bm{\hat{\xi}}_{11}, \bm{\hat{\xi}}_{21}, \ldots, \bm{\hat{\xi}}_{t K_t} \right)$ respectively, and each sub-vector of them, i.e., $\bm{\bar{\xi}}_{t'i} = \left( \bar{\xi}_{t'i}^1, \ldots, \bar{\xi}_{t'i}^{r_{t'i}} \right)$ and $\bm{\hat{\xi}}_{t'i} = \left( \hat{\xi}_{t'i}^1, \ldots, \hat{\xi}_{t'i}^{g_{t'i}} \right)$, can be derived by \eqref{eqn:MultiStageLiftedUncertainty1}. 
With $r_{ti}-1$ breakpoints placed inside the marginal support of each uncertain parameter $\xi_{ti}$, the two lifting operators defined in Subsection \ref{sec:LiftedUncertainty} can also be applied in the multistage case. We have $\bar{L}_t: \mathbb{R}^{K_t} \mapsto \mathbb{R}^{\overline{K}_t}$ with $\overline{K}_t = \sum_{i=1}^{K_{t}} r_{ti}$ and $\widehat{L}_t: \mathbb{R}^{{K}_t} \mapsto \mathbb{R}^{\widehat{K}_t}$ with $\widehat{K}_t = \sum_{i=1}^{K_{t}} g_{ti}$, where $g_{ti} = \max \left \{ 1, r_{ti}-1 \right \}$, for all $t=1, \ldots, T$. The lifted uncertain parameters $\bm{\bar{\xi}}_t = \left (\bm{\bar{\xi}}_{t1}, \ldots, \bm{\bar{\xi}}_{t K_t}  \right)$ and $\bm{\hat{\xi}}_t = \left (\bm{\hat{\xi}}_{t1}, \ldots, \bm{\hat{\xi}}_{t K_t}  \right)$ are then defined as follows:
\begin{subequations}
\label{eqn:MultiStageLiftedUncertainty1}
\begin{align}
    & \bar{\xi}_{ti}^j = \bar{L}_{ti}^j \left( \bm{e}_{ti} \xi_{ti}  \right) \quad \forall \, t=1, \ldots, T, \;  i =1,\ldots, K_{t}, \; j=1, \ldots, r_{ti} \\
    & \hat{\xi}_{ti}^j = \widehat{L}_{ti}^j  \left( \bm{e}_{ti} \xi_{ti}  \right) \quad \forall \, t=1, \ldots, T, \; i =1,\ldots, K_{t}, \; j=1, \ldots, g_{ti}.
\end{align}
\end{subequations}

Let $\bm{\tilde{\xi}}_t = \left( \bm{\xi}_t, \bm{\bar{\xi}}_t, \bm{\hat{\xi}}_t \right) \in \mathbb{R}^{\widetilde{K}_t}$ with $\widetilde{K}_t = K_t + \overline{K}_t + \widehat{K}_t$ and $\bm{\tilde{\xi}}^{[t]} = \left( \bm{\tilde{\xi}}_1, \ldots, \bm{\tilde{\xi}}_t \right)$. The lifted uncertainty set in stage $t$ is then
\begin{equation}
\label{eqn:MultiStageOpenLiftedUncertainty}
    \widetilde{\Xi}'^{[t]}(\bm{y}^{[t-1]}) = \left \{
    \bm{\tilde{\xi}}^{[t]} \in \mathbb{R}^{\widetilde{K}^{[t]}}: 
    \begin{array}{l}
         \bm{\xi}^{[t]} =  \left( \bm{\xi}_1, \bm{\xi}_2, \ldots, \bm{\xi}_t  \right) \in \Xi^{[t]}(\bm{y}^{[t-1]})  \\
         \bm{\bar{\xi}}_{t'} = \bar{L}_{t'}(\bm{\xi}_{t'}) \quad \forall \, t'=1,\ldots, t\\
         \bm{\hat{\xi}}_{t'} = \widehat{L}_{t'}(\bm{\xi}_{t'}) \quad \forall \, t'=1,\ldots, t
    \end{array}
    \right \},
\end{equation}
where $\widetilde{K}^{[t]} = \sum_{t'=1}^t \widetilde{K}_{t'}$.

Since $\widetilde{\Xi}'^{[t]}$ is an open set due to the discontinuity in $\bm{\hat{\xi}}_{t'}$, we also require the outer approximation of its closure, which we denote by $\widetilde{\Xi}^{[t]}$. For every $t=1,\ldots,T$ and $i=1, \ldots, K_t$, we define the following sets of vertices in the lifted space:
\begin{equation}
\label{eqn:VerticesSetMultiStage}
    \begin{array}{l}
         \widetilde{\mathcal{V}}_{ti} = \bigcup\limits_{j=1}^{r_{ti}} \widetilde{\mathcal{V}}_{ti}^j \\
         \widetilde{\mathcal{V}}_{ti}^j = \left \{
         \bm{\tilde{v}}_{ti} = \left( v_{ti}, \bm{\bar{v}}_{ti}, \bm{\hat{v}}_{ti} \right):
         \begin{array}{l}
              v_{ti} \in \mathcal{V}_{ti}^j  \\
              \bm{\bar{v}}_{ti} = \bar{L}_{ti} \left( \bm{e}_{ti} v_{ti} \right) \\
              \bm{\hat{v}}_{ti} = \lim\limits_{\xi_{ti} \rightarrow v_{ti}, \; \xi_{ti} \in \Xi_{ti}^j} \widehat{L}_{ti} \left( \bm{e}_{ti} \xi_{ti} \right)
         \end{array}
         \right \} \quad \forall \, j = 1, \dots, r_{ti},
    \end{array}
\end{equation}
where $\Xi_{ti}^j = \left\{\xi_{ti} \in \mathbb{R}: p_{ti}^{j-1} \leq \xi_{ti} \leq p_{ti}^j \right\}$ and $\mathcal{V}_{ti}^j = \left\{ p_{ti}^{j-1}, \, p_{ti}^j \right\}$, assuming that $p_{ti}^0 = \xi_{ti}^{\min}$ and $p_{ti}^{r_{ti}} = \xi_{ti}^{\max}$. This allows us to formulate the following convex hull representation of $\widetilde{\Xi}^{[t]}$:
\begin{equation}
\label{eqn:MultiStageLiftedUncertaintySet}
\begin{array}{l}
    \widetilde{\Xi}^{[t]}(\bm{y}^{[t-1]}) = \left \{
    \bm{\tilde{\xi}}^{[t]} \in \mathbb{R}^{\widetilde{K}^{[t]}}:
    \begin{array}{l}
         \sum\limits_{j=1}^{r_{t'i}} \sum\limits_{\bm{\tilde{v}}_{t'i} \in \widetilde{\mathcal{V}}_{t'i}^j} \lambda_{t'i}^j(\bm{\tilde{v}}_{t'i}) = 1 \quad \forall \, t' = 1, \ldots, t, \; i = 1, \ldots, K_{t'}  \\
         \xi_{t'i} = \sum\limits_{j=1}^{r_{t'i}} \sum\limits_{\bm{\tilde{v}}_{t'i} \in \widetilde{\mathcal{V}}_{t'i}^j}  \lambda_{t'i}^j(\bm{\tilde{v}}_{t'i}) v_{t'i} \quad \forall \, t' = 1, \ldots, t, \; i = 1, \ldots, K_{t'} \\
         \bm{\bar{\xi}}_{t'i} = \sum\limits_{j=1}^{r_{t'i}} \sum\limits_{\bm{\tilde{v}}_{t'i} \in \widetilde{\mathcal{V}}_{t'i}^j}  \lambda_{t'i}^j(\bm{\tilde{v}}_{t'i}) \bm{\bar{v}}_{t'i} \quad \forall \, t' = 1, \ldots, t, \; i = 1, \ldots, K_{t'} \\
         \bm{\hat{\xi}}_{t'i} = \sum\limits_{j=1}^{r_{t'i}} \sum\limits_{\bm{\tilde{v}}_{t'i} \in \widetilde{\mathcal{V}}_{t'i}^j}  \lambda_{t'i}^j(\bm{\tilde{v}}_{t'i}) \bm{\hat{v}}_{t'i} \quad \forall \, t' = 1, \ldots, t, \; i = 1, \ldots, K_{t'} \\
        \bm{W}^{[t]} \bm{\xi}^{[t]} \leq \sum\limits_{\hat{t}=1}^{t-1} \bm{U}_{\hat{t}}^{t} \bm{y}_{\hat{t}}(\bm{\tilde{\xi}}^{[\hat{t}]}) \\
         \lambda_{t'i}^j(\bm{\tilde{v}}_{t'i}) \in \mathbb{R}_+ \quad \forall \, t' = 1, \ldots, t, \; i = 1, \ldots, K_{t'}, \; j=1, \ldots, r_{t'i}
    \end{array}
    \right \}.
\end{array}
\end{equation}

\subsection{Decision Rule Approximation and Reformulation}
\label{sec:MultiStageReformulation}

Analogous to the decision rules \eqref{eqn:ProposedLinearRule} introduced in the two-stage problem, we apply the following decision rules to the recourse variables $\bm{x}_t$ and $\bm{y}_t$ for $t \geq 2$:
\begin{subequations}
\label{eqn:MultiStageDecisionRules}
\begin{align}
    & \bm{x}_{t} = \sum\limits_{t'=1}^t \sum\limits_{i=1}^{K_{t'}} \bm{\overline{X}}_{t'i}^t \bm{\bar{\xi}}_{t'i} + \bm{\widehat{X}}_{t'i}^t \bm{\hat{\xi}}_{t'i} \\
    & \bm{y}_t = \sum\limits_{t'=1}^t \sum\limits_{i=1}^{K_{t'}} \bm{\widehat{Y}}_{t'i}^t  \bm{\hat{\xi}}_{t'i}
\end{align}
\end{subequations}
where $\bm{\overline{X}}_{t'i}^t \in \mathbb{R}^{P_t \times r_{t'i}}$, $\bm{\widehat{X}}_{t'i}^t \in \mathbb{R}^{P_t \times g_{t'i}}$, and $\bm{\widehat{Y}}_{t'i}^t \in \{-1, \,0, \, 1\}^{Q_t \times g_{t'i}}$. Note that in the case of type-2 endogenous uncertainty, nonanticipativity has to be further enforced depending on which uncertain parameters materialize. This is encoded in the definition of the uncertainty set, where an unmaterialized uncertain parameter is forced to take the value zero. For linear decision rules like the ones in \eqref{eqn:MultiStageDecisionRules}, this is equivalent to forcing the decision rule coefficients corresponding to unmaterialized uncertain parameters to zero.

By substituting \eqref{eqn:LMIinFirstStage} and \eqref{eqn:MultiStageDecisionRules} into formulation \eqref{model:GeneralMultistageDeterministicModel} and replacing $\Xi^{[t]}(\bm{y}^{[t-1]})$ with $\widetilde{\Xi}^{[t]}(\bm{y}^{[t-1]})$, constraints \eqref{eqn:MultiStageDeterministicConstraint} become
\begin{equation}
\begin{aligned}
    \label{eqn:MultiStageSubstitutingConstraint}
    \sum\limits_{t'=1}^t \sum\limits_{i=1}^{K_{t'}} \bm{f}_{t'i}^t \xi_{t'i} +  \sum\limits_{t'' = 2} ^t \sum_{t'=1}^{t''} \sum\limits_{i=1}^{K_{t'}} \left[ \bm{\widetilde{A}}_{t''}^t \bm{\overline{X}}_{t'i}^{t''} \bm{\bar{\xi}}_{t'i} + \left( \bm{\widetilde{A}}_{t''}^t \bm{\widehat{X}}_{t'i}^{t''} + \bm{\widetilde{D}}_{t''}^t \bm{\widehat{Y}}_{t'i}^{t''}  \right) \bm{\hat{\xi}}_{t'i} \right] \leq \bm{0} \\
    \quad \forall \, t=2, \ldots, T,\; \bm{\tilde{\xi}}^{[t]} \in \widetilde{\Xi}^{[t]}(\bm{y}^{[t-1]})
\end{aligned}
\end{equation}
with $\bm{f}_{t'i}^t = \bm{A}_{t'i}^t \bm{x}_1 + \bm{D}_{t'i}^t \bm{y}_1 - \bm{b}_{t'i}^t$. The worst-case reformulation of \eqref{eqn:MultiStageSubstitutingConstraint} is then
\begin{equation}
\begin{aligned}
    \label{eqn:MultiStageWorst-CaseReformulation}
    \max_{\bm{\tilde{\xi}}^{[t]} \in \widetilde{\Xi}^{[t]}(\bm{y}^{[t-1]})} \;  \left \{ \sum\limits_{t'=1}^t \sum\limits_{i=1}^{K_{t'}} \left( \bm{f}_{t'i}^t \xi_{t'i} +  \sum\limits_{t'' = \max \{2, \, t' \}}^t \left[ \bm{\widetilde{A}}_{t''}^t \bm{\overline{X}}_{t'i}^{t''} \bm{\bar{\xi}}_{t'i} + \left( \bm{\widetilde{A}}_{t''}^t \bm{\widehat{X}}_{t'i}^{t''} + \bm{\widetilde{D}}_{t''}^t \bm{\widehat{Y}}_{t'i}^{t''}  \right) \bm{\hat{\xi}}_{t'i} \right] \right) \right \} \leq \bm{0} \\
    \forall \, t=2, \ldots, T,
\end{aligned}
\end{equation}
which in turn can be reformulated into the following set of constraints: \begin{subequations}
\label{eqn:MultiStageReformulatedConstraint}
    \begin{align}
        & \bm{\Phi}_t \bm{U}_1^t \bm{y}_1 + \sum\limits_{t'=1}^{t} \sum\limits_{i=1}^{K_{t'}} \bm{\delta}_{t'i}^t \leq \bm{0} \quad \forall \, t=2,\ldots, T  \label{eqn:MultiStageReformulatedConstraint-1}\\
        & \left(\bm{f}_{t'i}^t  - \bm{\Phi}_{t} \bm{w}_{t'i}^{[t]} \right) v_{t'i} - \bm{\delta}_{t'i}^t + \sum\limits_{t'' = \max\{2, \, t' \}}^t \left[ \bm{\widetilde{A}}_{t''}^t \bm{\overline{X}}_{t'i}^{t''} \bm{\bar{v}}_{t'i} + \left( \bm{\widetilde{A}}_{t''}^t \bm{\widehat{X}}_{t'i}^{t''} + \bm{\widetilde{D}}_{t''}^t \bm{\widehat{Y}}_{t'i}^{t''} + \bm{\Phi}_t \bm{U}_{t''}^t \bm{\widehat{Y}}_{t'i}^{t''}  \right) \bm{\hat{v}}_{t'i} \right] \leq \bm{0} \notag \\
        & \qquad\qquad\qquad\qquad\qquad\qquad \forall \, t=2,\ldots,T, \; t' = 1, \ldots, t, \; i = 1, \ldots, K_{t'}, \; \left( v_{t'i}, \, \bm{\bar{v}}_{t'i}, \, \bm{\hat{v}}_{t'i} \right) \in \widetilde{\mathcal{V}}_{t'i} \label{eqn:MultiStageReformulatedConstraint-2}\\
        & \bm{\Phi}_t \in \mathbb{R}_{+}^{N_t \times M_t}, \; \bm{\delta}_{t'i}^t \in \mathbb{R}^{N_t} \quad \forall \, t=2,\ldots,T, \; t'= 1, \ldots, t, \; i=1,\ldots, K_{t'} \label{eqn:MultiStageDualVariableDomain-1}
    \end{align}
\end{subequations}
where $\bm{U}_t^t = \bm{0}$ is introduced for notational convenience. The dual variables associated with the maximization problems in \eqref{eqn:MultiStageWorst-CaseReformulation} are denoted by $\bm{\Phi}_t$ and $\bm{\delta}_{t'i}^t$, and $\bm{w}_{t'i}^{[t]}$ is the column vector corresponding to $\xi_{t'i}$ in matrix $\bm{W}^{[t]}$. Constraints \eqref{eqn:MultiStageReformulatedConstraint-2} contain bilinear terms comprised of continuous and integer variables as $\bm{\widehat{Y}}_{t'i}^{t} \in \{-1, \, 0, \, 1\}^{Q_{t} \times g_{t'i}}$. To facilitate the linearization, we express $\bm{\widehat{Y}}_{t'i}^{t}$ as the difference between binary variables, i.e.
\begin{equation*}
    \bm{\widehat{Y}}_{t'i}^{t} = \bm{\dot{\widehat{Y}}}_{t'i}^{t} -\bm{\ddot{\widehat{Y}}}_{t'i}^{t}, \quad \bm{\dot{\widehat{Y}}}_{t'i}^{t} + \bm{\ddot{\widehat{Y}}}_{t'i}^{t} \leq \bm{1}
\end{equation*}
where $\bm{\dot{\widehat{Y}}}_{t'i}^{t}, \, \bm{\ddot{\widehat{Y}}}_{t'i}^{t} \in \{0,1\}^{Q_{t} \times g_{t'i}}$, and the latter inequalities can be added to eliminate symmetry. Constraints \eqref{eqn:MultiStageReformulatedConstraint-2} then become
\begin{equation}
    \label{eqn:MultiStageTransformedReformulatedConstraint-2}
    \begin{aligned}
    & \left(\bm{f}_{t'i}^t  - \bm{\Phi}_{t} \bm{w}_{t'i}^{[t]} \right) v_{t'i} - \bm{\delta}_{t'i}^t \\
    & + \sum\limits_{t'' = \max\{2, t' \}}^t \left( \bm{\widetilde{A}}_{t''}^t \bm{\overline{X}}_{t'i}^{t''} \bm{\bar{v}}_{t'i} + 
    % \sum\limits_{t'' = \max\{2, t' \}}^t 
    \left[ \bm{\widetilde{A}}_{t''}^t \bm{\widehat{X}}_{t'i}^{t''} + \left (\bm{\widetilde{D}}_{t''}^t  + \bm{\Phi}_t \bm{U}_{t''}^t \right) \left( \bm{\dot{\widehat{Y}}}_{t'i}^{t''} -\bm{\ddot{\widehat{Y}}}_{t'i}^{t''} \right)  \right] \bm{\hat{v}}_{t'i} \right) \leq \bm{0}  \\
    & \qquad\qquad\qquad\qquad \forall \, t=2,\ldots,T, \; t' = 1, \ldots, t, \; i = 1, \ldots, K_{t'}, \; \left( v_{t'i}, \, \bm{\bar{v}}_{t'i}, \, \bm{\hat{v}}_{t'i} \right) \in \widetilde{\mathcal{V}}_{t'i}
    \end{aligned}
\end{equation}

The integrality constraints on $\bm{y}_t$ in \eqref{eqn:MultiStageDomainConstraint} can be relaxed, i.e. $\bm{0} \leq \bm{y}_t(\bm{\xi}^{[t]}) \leq \bm{e}$, and reformulated as follows:
\begin{equation}
\label{eqn:MultiStageDomainReformulation}
    \begin{aligned}
     & \left .
    \begin{array}{c}
        \sum\limits_{t'=1}^t \sum\limits_{i=1}^{K_{t'}} \bm{\underline{\pi}}_{t'i}^t - \bm{\underline{\Theta}}^t \bm{U}_1^t \bm{y}_1 \geq  \bm{0}   \\
        \left( \bm{\dot{\widehat{Y}}}_{t'i}^{t} -\bm{\ddot{\widehat{Y}}}_{t'i}^{t} \right) \bm{\hat{v}}_{t'i} - \bm{\underline{\pi}}_{t'i}^t + \bm{\underline{\Theta}}^t \bm{w}_{t'i}^{[t]} v_{t'i} - \sum\limits_{t'' = \max\{2, \, t' \}}^t \bm{\underline{\Theta}}^t \bm{U}_{t''}^t \left(\bm{\dot{\widehat{Y}}}_{t'i}^{t''} -\bm{\ddot{\widehat{Y}}}_{t'i}^{t''} \right) \bm{\hat{v}}_{t'i}  \geq  \bm{0}  \\
        \qquad \qquad \qquad \qquad \qquad \qquad \quad   \forall \, t' = 1, \ldots, t, \; i= 1, \ldots, K_{t'}, \; \left( v_{t'i}, \, \bm{\bar{v}}_{t'i}, \, \bm{\hat{v}}_{t'i} \right) \in \widetilde{\mathcal{V}}_{t'i} \\
        \bm{\underline{\Theta}}^t \in \mathbb{R}_+^{Q_t \times M_t}, \; \bm{\underline{\pi}}_{t'i}^t \in \mathbb{R}^{Q_t} \quad \forall \, t' = 1, \ldots, t, \; i= 1, \ldots, K_{t'}
    \end{array}
    \right \}  \forall \, t=2, \ldots, T  \\
    & \left . \begin{array}{c}
        \sum\limits_{t'=1}^t \sum\limits_{i=1}^{K_{t'}} \bm{\overline{\pi}}_{t'i}^t + \bm{\overline{\Theta}}^t \bm{U}_1^t \bm{y}_1 \leq \bm{e}  \\
        \left( \bm{\dot{\widehat{Y}}}_{t'i}^{t} -\bm{\ddot{\widehat{Y}}}_{t'i}^{t} \right) \bm{\hat{v}}_{t'i} - \bm{\overline{\pi}}_{t'i}^t - \bm{\overline{\Theta}}^t \bm{w}_{t'i}^{[t]} v_{t'i} + \sum\limits_{t'' = \max\{2, \, t' \}}^t \bm{\overline{\Theta}}^t \bm{U}_{t''}^t \left(\bm{\dot{\widehat{Y}}}_{t'i}^{t''} -\bm{\ddot{\widehat{Y}}}_{t'i}^{t''} \right) \bm{\hat{v}}_{t'i}  \leq \bm{0} \\
        \qquad \qquad \qquad \qquad \qquad \qquad \quad \forall \, t' = 1, \ldots, t, \; i= 1, \ldots, K_{t'}, \; \left( v_{t'i}, \, \bm{\bar{v}}_{t'i}, \, \bm{\hat{v}}_{t'i} \right) \in \widetilde{\mathcal{V}}_{t'i}\\
        \bm{\overline{\Theta}}^t \in \mathbb{R}_+^{Q_t \times M_t}, \;\bm{\overline{\pi}}_{t'i}^t  \in \mathbb{R}^{Q_t} \quad \forall \, t' = 1, \ldots, t, \; i= 1, \ldots, K_{t'}
    \end{array}
    \right \} \forall \, t=2, \ldots, T.
    \end{aligned}
\end{equation}
%where $\bm{\underline{\pi}}_{t'i}^t$, $\bm{\underline{\Theta}}^t$, $\bm{\overline{\pi}}_{t'i}^t  $ and $\bm{\overline{\Theta}}^t$ are dual variables. 

Finally, we arrive at the following reformulation of the multistage problem:
\begin{subequations}
\label{eqn:MultiStageReformulation}
\begin{align}
    \minimize \quad & \bm{e}_1^{\top} \bm{x}_1 \\
    \st \quad & \bm{A}^1 \bm{x}_1 + \bm{D}^1 \bm{y}_1 \leq \bm{b}^1 \\
    & \bm{f}_{t'i}^t = \bm{A}_{t'i}^t \bm{x}_1 + \bm{D}_{t'i}^t \bm{y}_1 - \bm{b}_{t'i}^t \quad \forall \, t= 2, \ldots, T, \, t' = 1, \ldots, t, \; i=1, \ldots, K_{t'} \\
    & \eqref{eqn:MultiStageReformulatedConstraint-1}, \; \eqref{eqn:MultiStageDualVariableDomain-1}, \; \eqref{eqn:MultiStageTransformedReformulatedConstraint-2}, \; \eqref{eqn:MultiStageDomainReformulation} \\
    & \bm{x}_1 \in \mathbb{R}^{P_1}, \; \bm{y}_1 \in \{ 0, 1 \}^{Q_1} \\
    & \bm{\overline{X}}_{t'i}^t \in \mathbb{R}^{P_t \times r_{t'i}}, \; \bm{\widehat{X}}_{t'i}^t \in \mathbb{R}^{P_t \times g_{t'i}}, \;  \bm{\dot{\widehat{Y}}}_{t'i}^{t}, \, \bm{\ddot{\widehat{Y}}}_{t'i}^{t} \in \{0,1\}^{Q_{t} \times g_{t'i}} \notag \\
    & \quad\quad\quad\quad\quad\quad\quad\quad\quad\quad \forall \, t=2, \ldots, T,\; t' = 1, \ldots, t, \; i=1, \ldots, K_{t'} .
\end{align}
\end{subequations}
Every bilinear term in \eqref{eqn:MultiStageReformulation} is composed of a continuous and a binary variable; hence, an MILP formulation is readily obtained after exact linearization of the bilinear terms.

Note that in the decision rules given by \eqref{eqn:MultiStageDecisionRules}, a recourse variable in stage $t$ is a function of all uncertain parameters observed up to stage $t$, i.e. the decision rules make use of all available information. In practice, however, this may not be the best choice as the model size and hence the computational performance strongly depend on the number of parameters involved in the decision rules. It has been observed in several multistage applications \citep{Zhang2016e, Lappas2016} that the optimal decision rules usually only depend on a small subset of uncertain parameters; hence, a common strategy to reduce computation time is to restrict the decision rules to depend on a smaller set of uncertain parameters. An intuitive choice is to let a recourse variable in stage $t$ only depend on the uncertain parameters observed in the previous $\Delta t$ and the current stages. The required change in the reformulation is shown in Appendix \ref{Appendix:DeltaTApproximation}.

\section{Computational Case Studies}
\label{sec:CaseStudies}

In this section, we apply the proposed approach to a two-stage design problem and to a multistage production planning problem. All model instances were implemented in Julia v1.2.0 using the modeling environment JuMP v0.20.0 \citep{Lubin2015} and solved to 1\% optimality using Gurobi v8.1.1 on a Intel Core i7-8700 CPU at 3.20 GHz machine with 8 GB RAM.

\subsection{Design for Flexible Production}

We consider the design of a production system that manufactures a single product for which the required production amount can vary across a wide range. Such production flexibility is especially important in systems with little or no product inventory capacity, e.g. in electricity generation or if the product is highly volatile. The production system can consist of a set of production units that all produce the same product but differ in capacity and cost. The production cost of a unit is assumed to be an affine function of the production amount. Moreover, while the minimum production amount, $c^{\min}$, is known, the maximum production amount (i.e. the production capacity) is uncertain and is only known after the unit is built. However, we do know that the capacity will be between $c^{\max} - \hat{c}^{\max}$ and $c^{\max}$.

%\begin{figure}[ht]\centering
%\includegraphics[width=3.5in]{Design1.pdf}
%\caption{For a given production unit, the production cost is an affine function of the production amount. The production capacity is uncertain until the unit is built, but will lie between $c^{\max} - \hat{c}^{\max}$ and $c^{\max}$.}
%\label{fig:Design1}
%\end{figure}

Given a set of production units $\mathcal{I}$, the objective is to decide which subset of units to build such that all product demand $d$ within a range $\left [d^{\min}, d^{\max} \right]$ can be met exactly and the worst-case total cost is minimized. We can formulate the problem as the following two-stage robust optimization problem:
\begin{subequations}
\label{eqn:Design}
\begin{align}
\minimize \quad & \sum_{i \in \mathcal{I}} \alpha_i z_i + \max_{\bm{\hat{c}} \in \widehat{\mathcal{C}}(\bm{z}), d \in \mathcal{D}} \sum_{i \in \mathcal{I}} \beta_i \tilde{y}_i(\bm{\hat{c}},d) + \gamma_i \tilde{x}_i(\bm{\hat{c}},d) \\
\st \quad & z_i \in \{0,1\} \quad \forall \, i \in \mathcal{I} \\
& \sum_{i \in \mathcal{I}} \tilde{x}_i(\bm{\hat{c}},d) = d \quad \forall \, \bm{\hat{c}} \in \widehat{\mathcal{C}}(\bm{z}), \, d \in \mathcal{D} \label{eqn:DesignEquality} \\
& \hspace{-8pt} \left.
\begin{array}{l}
  \tilde{y}_i(\bm{\hat{c}},d) \leq z_i \\[4pt]
  \tilde{x}_i(\bm{\hat{c}},d) \geq c^{\min}_i \tilde{y}_i(\bm{\hat{c}},d) \\[4pt]
  \tilde{x}_i(\bm{\hat{c}},d) \leq c^{\max}_i \tilde{y}_i(\bm{\hat{c}},d) \\[4pt]
  \tilde{x}_i(\bm{\hat{c}},d) \leq c^{\max}_i - \hat{c}_i \\[4pt]
  \tilde{x}_i(\bm{\hat{c}},d) \in \mathbb{R}_+ \\[4pt]
  \tilde{y}_i(\bm{\hat{c}},d) \in \{0,1\}
\end{array}
\right\rbrace \quad \forall \, i \in \mathcal{I}, \, \bm{\hat{c}} \in \widehat{\mathcal{C}}(\bm{z}), \, d \in \mathcal{D}
\end{align}
\end{subequations}
with the uncertainty sets $\widehat{\mathcal{C}}(\bm{z})$ and $\mathcal{D}$ defined as follows:
\begin{equation*}
\widehat{\mathcal{C}}(\bm{z}) = \left\lbrace \bm{\hat{c}} \in \mathbb{R}_+^{|\mathcal{J}|}: \bm{\hat{c}} \leq \bm{\hat{c}}^{\max} \circ \bm{z} \right\rbrace
\end{equation*}
\begin{equation*}
\mathcal{D} = \left\lbrace d \in \mathbb{R}_+: d^{\min} \leq d \leq d^{\max} \right\rbrace.
\end{equation*}

In problem \eqref{eqn:Design}, $\bm{z}$ are the first-stage design variables while $\bm{\tilde{x}}$ and $\bm{\tilde{y}}$ are the second-stage operational variables. Production unit $i$ is built if $z_i = 1$; $\tilde{y}_i = 1$ if unit $i$ is used to manufacture the product, and $\tilde{x}_i$ denotes the corresponding production amount. In the objective function, $\alpha_i$, $\beta_i$, and $\gamma_i$ denote the capital, fixed production, and variable production costs for unit $i$, respectively.

\subsubsection{The Benefit of Discrete Recourse}

\begin{figure}[ht]\centering
\includegraphics[width=5in]{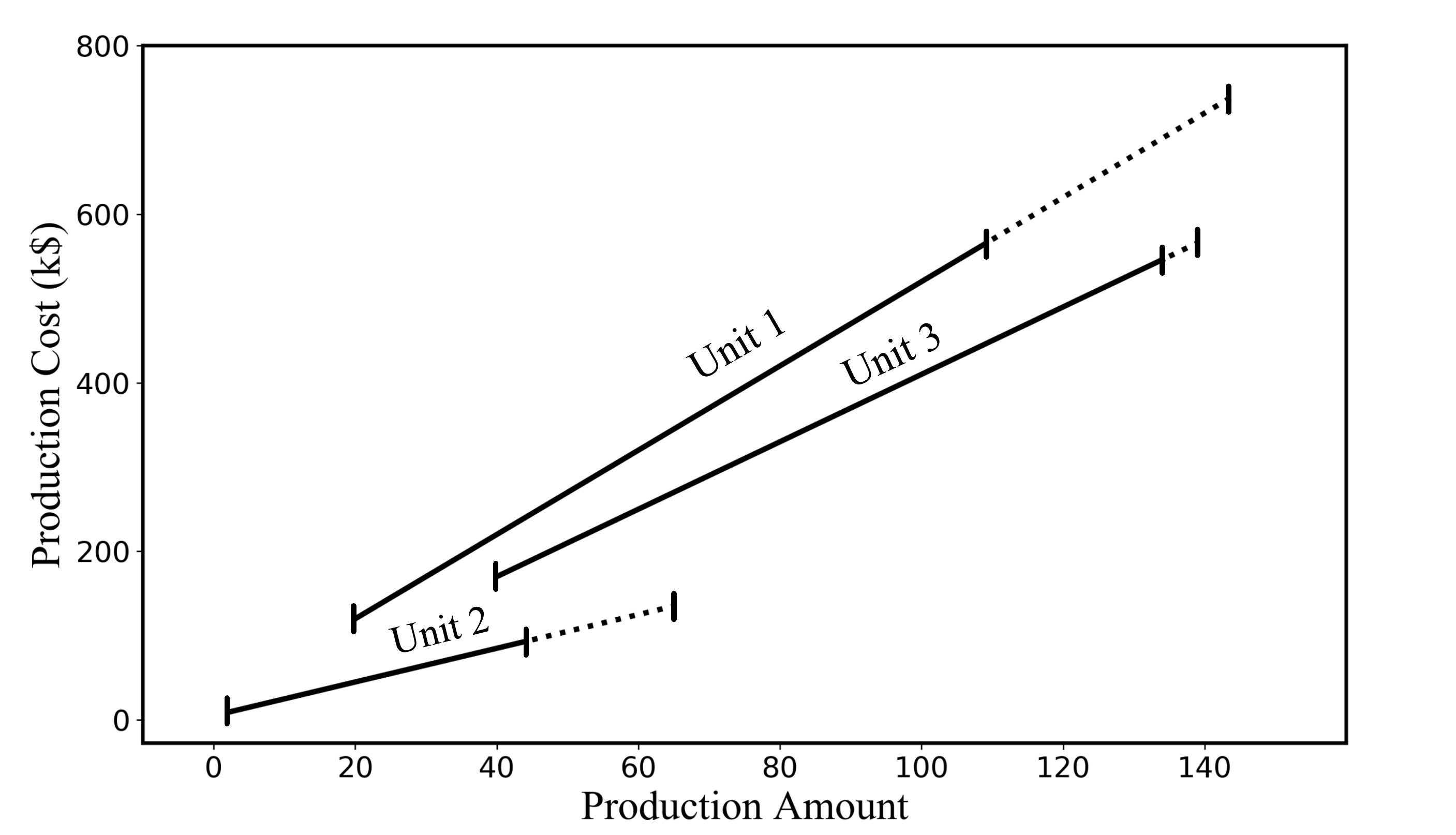}
\caption{Capacity ranges and production cost functions of Units 1--3.}
\label{fig:Design2}
\end{figure}

Demand $d$ is modeled as an exogenous uncertain parameter. Notice that because of the equality constraints \eqref{eqn:DesignEquality}, the problem will be infeasible if we apply static robust optimization, i.e. if we assume that $\bm{\tilde{x}}$ and $\bm{\tilde{y}}$ are not adjustable. In the following, we further demonstrate the importance of discrete recourse with an illustrative example, for which all data are provided in Appendix \ref{apx:Data}. It involves three alternative production units with different capacity ranges and production cost functions as shown in Figure \ref{fig:Design2}, where the dashed line segments indicate the regions of uncertainty. Consider two cases with different demand ranges: Case A with $d^{\min} = \min_{i \in \mathcal{I}} c_i^{\min}$ and $d^{\max} = \sum_{i \in \mathcal{I}} (c_i^{\max} - \hat{c}_i^{\max})$, and Case B with $d^{\min} = c^{\min}_1$ and $d^{\max} = c_1^{\max} - \hat{c}_1^{\max}$.

% \begin{figure}[ht]\centering
%     \subfloat[]{
%         \label{fig: two-stage case 1-1}
%         \fbox{\includegraphics[width = 3in]{Figures/TwoStage1-1.pdf}}}
%     \subfloat[]{
%         \label{fig: two-stage case 1-2}
%         \fbox{\includegraphics[width = 3in]{Figures/TwoStage1-2.pdf}}}
% \caption{Capacity ranges and production cost functions of Units 1--3.}
% \label{fig:Design2}
% \end{figure} 

In Case A, the maximal possible demand $d^{\max}$ can be fulfilled by installing and operating all three units, while only Unit 2, which has the smallest capacity, should be operated to meet the minimal demand $d^{\min}$. The problem is infeasible if only continuous recourse is considered since the decision of which units to operate cannot be adjusted depending on the realization of the demand. However, if the binary variables $\bm{\tilde{y}}$ are adjustable, any demand within the given range can be met exactly.

In Case B, the demand range is such that it can be covered by Unit 1 but neither Unit 2 nor Unit 3 alone. However, with both Units 2 and 3 installed and switching between these two units depending on the realization of the demand, meeting demand over the entire range is feasible. Moreover, the production costs for Units 2 and 3 are lower than for Unit 1. Hence, assuming the capital costs are the same, selecting Units 2 and 3 is a better solution than selecting Unit 1. Yet, this solution is only feasible if we have both binary and continuous recourse, i.e. if $\bm{\tilde{x}}$ and $\bm{\tilde{y}}$ are both adjustable. If only $\bm{\tilde{x}}$ are adjustable, the only robust feasible solution is to choose Unit 1.

Table \ref{tab:table 1 of two-stage case} shows the optimal (worst-case) costs for Cases A and B, each solved once allowing only continuous recourse and then with both continuous and binary recourse. The decision rules are constructed using three equidistantly generated breakpoints for each uncertain parameter in both Case A and B. One can see that Case A is only feasible if we allow both continuous and binary recourse. In Case B, considering binary in addition to continuous recourse significantly improves the optimal value while ensuring the same level of robustness.

\begin{table}[h!]\centering
\setlength\tabcolsep{10pt}	
	\caption{Optimal (worst-case) costs for Cases A and B, solved with only continuous and with both continuous and binary recourse.}
	\begin{tabular}{ccc}
	\toprule
	\textbf{Case} & \textbf{Recourse} & \textbf{Optimal value (k\$)} \\ \midrule
	\multirow{2}{*}{A}  & continuous only & infeasible \\
	                    & continuous \& binary & 1,415 \\ \midrule
    \multirow{2}{*}{B}  & continuous only & 670 \\
	                    & continuous \& binary & 465 \\
	\bottomrule
	\end{tabular}
	\label{tab:table 1 of two-stage case}
\end{table}
 
% {\color{Maroon} Show example data in appendix. Solve problem with continuous recourse only and then with both continuous and integer recourse.}

% \begin{equation}
% \begin{array}l
%     z_1=0, \; z_2=1, \; z_3=1 \\
%     \tilde{y}_1=0, \; \tilde{x}_1 = 0 \\
%     s
% \end{array}
% \end{equation}
 
\subsubsection{The Impact of Endogenous Uncertainty}

In this problem, we model $\bm{\hat{c}}$ as endogenous uncertain parameters. However, given the uncertainty set $\widehat{\mathcal{C}}(\bm{z})$, it is actually not necessary to do so. The main reason is that since the uncertain parameters are independent, $\hat{c}_i = \hat{c}^{\max}_i$ is always the worst case if $z_i = 1$. In addition, the value of $\hat{c}_i$ does not affect the problem if $z_i = 0$. Hence, we can simply replace $\bm{\hat{c}}$ in \eqref{eqn:Design} with $\bm{\hat{c}}^{\max}$. The situation is different if the uncertain parameters are correlated. For example, consider the following budget-based uncertainty set for $\bm{\hat{c}}$:
\begin{equation*}
\widehat{\mathcal{C}}^1 = \left\lbrace \bm{\hat{c}} \in \mathbb{R}_+^{|\mathcal{I}|}: \sum_{i \in \mathcal{I}} \hat{c}_i \leq \tau \sum_{i \in \mathcal{I}} \hat{c}_i^{\max}, \; \bm{\hat{c}} \leq \bm{\hat{c}}^{\max}\right\rbrace,
\end{equation*}
where the total deviation of $\bm{\hat{c}}$ from zero across all possible units is bounded from above by $\tau \sum_{i \in \mathcal{I}} \hat{c}_i^{\max}$. In this case, for $\tau < 1$, the worst-case value of $\hat{c}_i$ may not be $\hat{c}^{\max}_i$ and cannot be easily determined a priori. 

The uncertain parameter $\hat{c}_i$ is endogenous because it only materializes if unit $i$ is built. Otherwise, the parameter is physically meaningless and should therefore be irrelevant for the problem; as a result, the budget uncertainty set should change accordingly. This endogenous nature of the uncertainty is not capture in $\widehat{\mathcal{C}}^1$. A more appropriate decision-dependent uncertainty set is
\begin{equation*}
    \widehat{\mathcal{C}}^2 (\bm{z}) = \left\lbrace \bm{\hat{c}} \in \mathbb{R}_+^{|\mathcal{I}|}: \sum_{i \in \mathcal{I}} \hat{c}_i \leq \tau \sum_{i \in \mathcal{I}} \hat{c}_i^{\max} z_i, \; \bm{\hat{c}} \leq \bm{\hat{c}}^{\max} \circ \bm{z}\right\rbrace,
\end{equation*}
where $\hat{c}_i$ is fixed to zero if $z_i=0$ and the budget only considers materialized uncertain capacities. One can see that, assuming $\bm{\hat{c}}^{\max} > \bm{0}$, $\widehat{\mathcal{C}}^2 (\bm{z}) \subset \widehat{\mathcal{C}}^1$ for any $\bm{z} \neq \bm{e}$; hence, using $\widehat{\mathcal{C}}^1$ as the uncertainty set is expected to lead to overly conservative solutions.

\begin{figure}[ht]\centering
\includegraphics[width=5in]{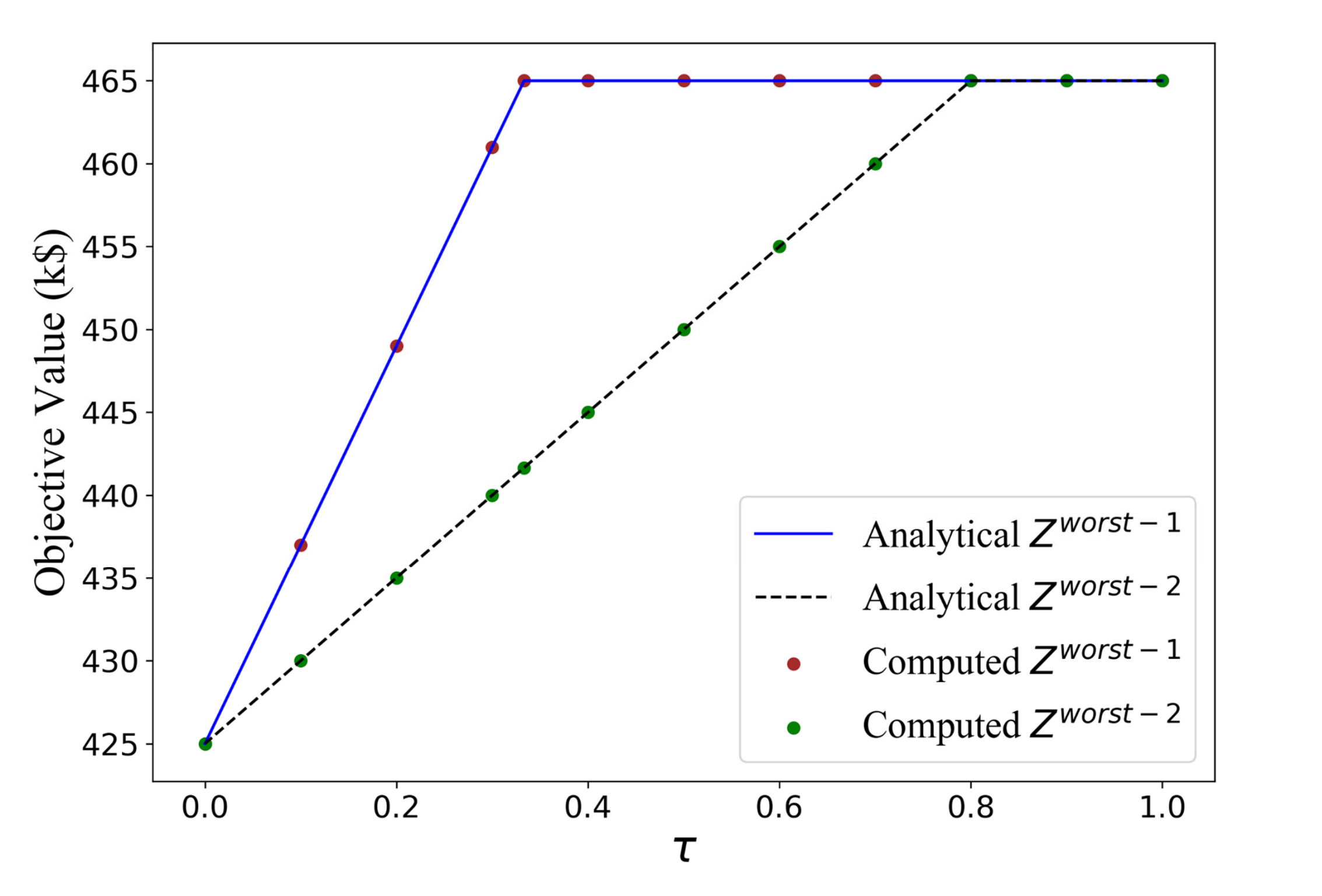}
\caption{Dependence of the optimal value in Case B on the choice of uncertainty set and the budget parameter $\tau$. Results also show that the proposed decision rules can achieve the optimal solutions with appropriate breakpoints.}
\label{fig:TwoStage ex2}
\end{figure}

Consider Case B for which we can obtain the following analytical optimal solution:
\begin{equation*}
\begin{aligned}
    & z_1 = 0, \; z_2 = z_3 = 1, \; \tilde{y}_1 = 0, \; \tilde{x}_1 = 0 \\
    & \left \{
        \begin{array}{ll}
             \tilde{y}_2 = 1, \; \tilde{x}_2 = d, \; \tilde{y}_3 = \tilde{x}_3 = 0, \quad & \text{if} \; d^{\min} \leq d \leq c_2^{\max} - \hat{c}_2 \\
             \tilde{y}_2 = \tilde{y}_3 = 1, \; \tilde{x}_2 = d - c_3^{\min}, \; \tilde{x}_3 = c_3^{\min}, \quad & \text{if} \; c_2^{\max} - \hat{c}_2 < d \leq c_3^{\min} + c_2^{\max} - \hat{c}_2 \\
             \tilde{y}_2 = \tilde{y}_3 = 1, \; \tilde{x}_2 = c_2^{\max} - \hat{c}_2, \; \tilde{x}_3 = d - c_2^{\max} + \hat{c}_2, \quad & \text{if} \; c_3^{\min} + c_2^{\max} - \hat{c}_2 < d \leq d^{\max}.
        \end{array}
    \right.
\end{aligned}
\end{equation*}
Then, we can determine the worst case, which depends on the choice of uncertainty set:
\begin{equation*}
\begin{aligned}
    & \text{For } \bm{\hat{c}} \in \widehat{\mathcal{C}}^1: \quad && d^{\text{worst-1}} = d^{\max}, \; \hat{c}_2^{\text{worst-1}} = \min \left \{\hat{c}_2^{\max}, \;  \tau \left (  \hat{c}_1^{\max} + \hat{c}_2^{\max} + \hat{c}_3^{\max}\right ) \right \}. \\
    & \text{For } \bm{\hat{c}} \in \widehat{\mathcal{C}}^2(\bm{z}): \quad && d^{\text{worst-2}} = d^{\max}, \; \hat{c}_2^{\text{worst-2}} =\min \left \{\hat{c}_2^{\max}, \; \tau \left ( \hat{c}_2^{\max} + \hat{c}_3^{\max} \right) \right \}.
\end{aligned}
\end{equation*}
The optimal values (in k\$), i.e. the minimum worst-case costs, for these two cases can then be computed as follows:
\begin{equation*}
\begin{aligned}
    & Z^{\text{worst-1}} = -15 + 4d^{\max} + 2 \hat{c}_2^{\text{worst-1}}. \\
    & Z^{\text{worst-2}} = -15 + 4d^{\max} + 2 \hat{c}_2^{\text{worst-2}}.
\end{aligned}
\end{equation*}
Evidently, we have $Z^{\text{worst-1}} \geq Z^{\text{worst-2}}$, which shows that the decision-dependent uncertainty set $\widehat{\mathcal{C}}^2(\bm{z})$ is less conservative than the fixed $\widehat{\mathcal{C}}^1$. We can further see that $Z^{\text{worst-1}}$ and $Z^{\text{worst-2}}$ both depend on $\tau$, which is depicted in Figure \ref{fig:TwoStage ex2}. Here, we also compare the analytical solutions with solutions obtained from solving the two-stage robust optimization problem with the proposed decision rules. The breakpoints for demand $d$ are chosen to be $c_2^{\max} - \hat{c}_2^{\text{worst}}$ and $c_3^{\min} + c_2^{\max} - \hat{c}_2^{\text{worst}}$ with $\hat{c}_2^{\text{worst}} = \hat{c}_2^{\text{worst-1}}$ or $\hat{c}_2^{\text{worst}} = \hat{c}_2^{\text{worst-2}}$ depending on the choice of uncertainty set. By doing so, as shown in Figure \ref{fig:TwoStage ex2}, we can recover the analytical optimal solutions.

\subsubsection{On the Selection of Breakpoints}

The quality of the proposed decision rules strongly depends on the choice of breakpoints. In Case B, the optimal breakpoints could be determined a priori; however, this is not generally true in more complex instances. In practice, we have to apply some heuristic to generate the breakpoints. The most intuitive one is to simply choose the number of breakpoints for each uncertain parameter and place them equidistantly inside its marginal support. However, it is recommended to utilize problem-specific features to design improved breakpoint generation procedures. For example, in this problem, according to constraint \eqref{eqn:DesignEquality}, the demand $d$ has to be equal to the sum of all built units' production amounts. This insight motivates a tailored method that uses all $c_i^{\min}$ and $c_i^{\max}$ that are within the range $[d^{\min}, d^{\max}]$ as breakpoints for $d$. More generally, let ${p}_i^j$ be a breakpoint for $\hat{c}_i$, choose all $c_i^{\min} \in [d^{\min}, d^{\max}]$ and $c_i^{\max} - {p}_i^j \in [d^{\min}, d^{\max}]$ for $i \in \mathcal{I}$ and $j=0, \ldots, {r}_i$ to be breakpoints for $d$.

For a randomly generated case with eight alternative production units (data provided in Appendix \ref{apx:Data}), we compare the two heuristic breakpoint generation methods described above. In the case of equidistant construction of breakpoints, we apply the same number of breakpoints to each of the nine uncertain parameters and examine the performance as we increase the number of breakpoints. When using the tailored method, we only apply breakpoints to $d$. The computational results are obtained with $d^{\min} = 43.1$, $d^{\max} = 406.5$ and $\tau = 0.5$, as shown in Table \ref{table:RandomTest}. One can see that the problem is infeasible if no breakpoints are used. Similarly, it is infeasible if one breakpoint is placed at the center of the marginal support of every uncertain parameter. As the number of breakpoints increases, equidistant generation of breakpoints leads to improved solutions, albeit at higher computational cost since the model size grows with the number of breakpoints. Note, however, that the optimal value does not improve from 27 to 36 breakpoints. The same optimal value is achieved with the 15 breakpoints generated using the tailored method, where the problem was solved in 12 seconds, which is in contrast to the 172 seconds required to solve the instance with a total of 27 equidistantly placed breakpoints. 

\begin{table}[h!]\centering
\setlength\tabcolsep{5pt}	
	\caption{Computational results for instances with breakpoints generated using different heuristics.}
	\begin{tabular}{ccccccc}
	\toprule
	\textbf{Method} & \makecell{\textbf{Total \# of} \\ \textbf{breakpoints}} & \makecell{\textbf{Optimal} \\ \textbf{value (\$)}} & \makecell{\textbf{Solution} \\ \textbf{time (s)}} & \makecell{\textbf{\# of} \\ \textbf{constraints}} & \makecell{\textbf{\# of continuous} \\ \textbf{variables}} & \makecell{\textbf{\# of integer} \\ \textbf{variables}} \\ \midrule
	& 0 & infeasible & n/a & 38,796 & 11,611 & 16 \\ \midrule
	\multirow{4}{*}{Equidistant}  
	& 9 & infeasible & n/a & 39,894 & 11,755 & 88 \\
	& 18 & 1,556,563 & 77.60 & 40,992 & 11,899 & 160 \\
	& 27 & 1,525,679 & 172.08 & 42,090 & 12,043 & 232 \\
	& 36 & 1,525,679 & 205.93 & 43,188 & 12,187 & 304 \\	\midrule
    Tailored & 15 & 1,525,679 & 12.32 &	40,626 & 11,851 & 136 \\
	\bottomrule
	\end{tabular}
	\label{table:RandomTest}
\end{table}

\subsection{Multiperiod Production Planning}

In the second case study, we consider a multiperiod production planning problem with exogenous uncertain demands and endogenous uncertain production capacities. Endogenous, especially type-2 endogenous, uncertainty is prevalent in planning and scheduling applications as many task-related uncertainties, such as production capacity, yield, and processing time, only materialize if one decides to perform the task \citep{Goel2004, Colvin2008, Lappas2016}.

\begin{figure}[h!]\centering
\includegraphics[width=4.5in]{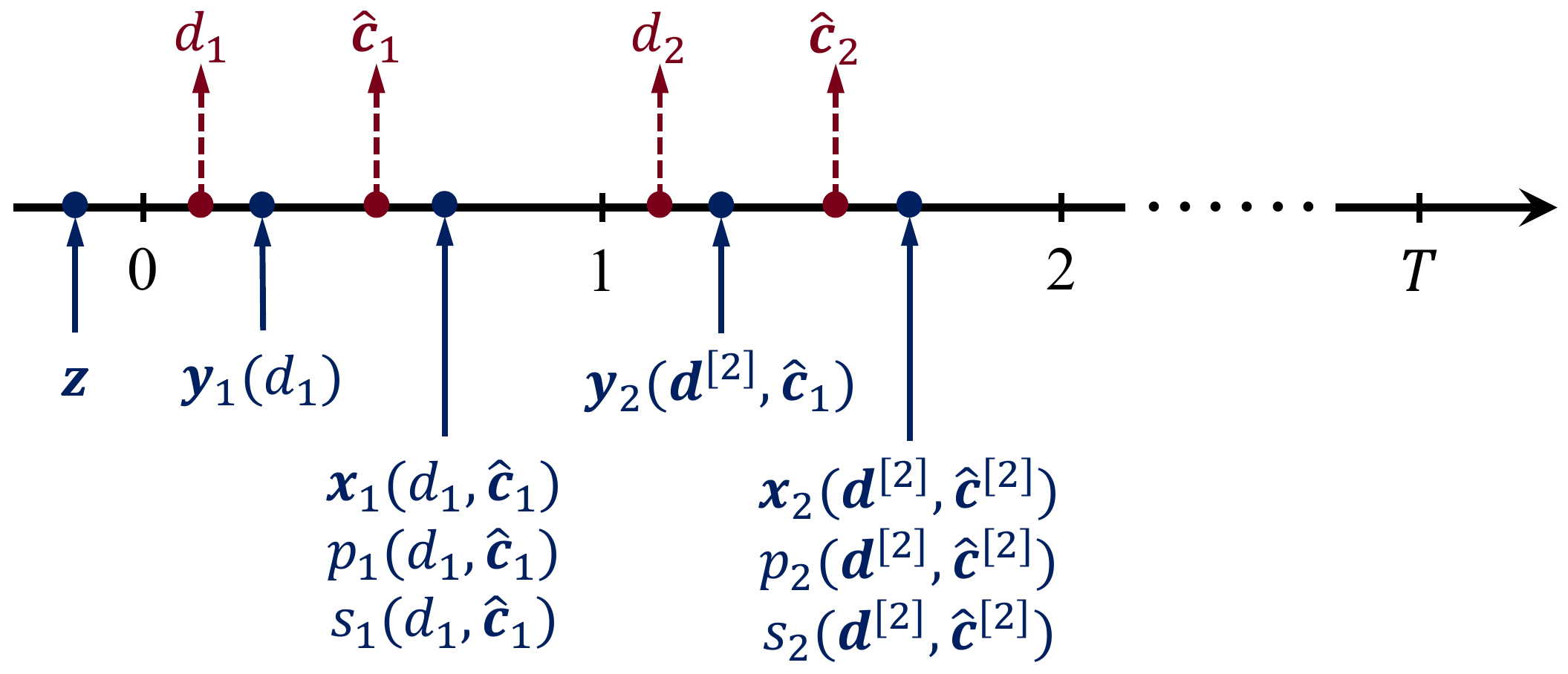}
\caption{The sequential decision-making process considered in the multistage formulation. Realizations of uncertain parameters are shown above the horizontal time axis, while decision variables are shown below the axis.}
\label{fig:MultistageStageVsPeriods}
\end{figure}

The multistage sequential decision-making process is depicted in Figure \ref{fig:MultistageStageVsPeriods}, where we apply the convention that a time period $t$ starts at time point $t-1$ and ends at time point $t$. Before the start of the planning horizon, which is given by the set of time periods $\mathcal{T} := \{1, \ldots, T \}$, we have to decide whether each unit $i \in \mathcal{I} := \{1, \ldots, I\}$ should be upgraded such that its capacity is increased or the uncertainty associated with the capacity is reduced. This first-stage binary decision is denoted by $z_i$ and is associated with a fixed cost $\gamma_i$. We then observe the demand in time period 1, $d_1$, and decide on which units to run; hence, the binary variable $y_{1i}$, which is 1 if and only if unit $i$ operates in time period 1, depends on the realization of $d_1$. The production capacity of unit $i$ in time period $t$ is $c_i^{\max} - \hat{c}_{ti}$, where $\hat{c}_{ti}$ is an uncertain parameter. The uncertainty in the capacity of a production unit only materializes if the unit is turned on; hence, $\bm{\hat{c}}_1$ are only observed after $\bm{y}_1$ are set. Once $\bm{\hat{c}}_1$ are observed, the production amounts $\bm{x}_1$, the purchasing amount $p_1$, and the resulting inventory level $s_1$ are determined; hence, these decisions depend on the realization of $d_1$ and $\bm{\hat{c}}_1$. As indicated in Figure \ref{fig:MultistageStageVsPeriods}, this sequential decision-making process is carried out until the end of the planning horizon. As a result, we have a multistage problem with $2T+1$ stages.

The multistage robust production planning problem is formulated as follows:
\begin{equation}
\label{Model:MultistageCase}
\begin{aligned}
\minimize \quad & \sum\limits_{i \in \mathcal{I}} \gamma_i z_i + \max_{\bm{d}^{[T]} \in \mathcal{D}^{[T]}, \, \bm{\hat{c}}^{[T]} \in \widehat{\mathcal{C}}^{[T]}(\bm{z}, \, \bm{y}^{[T]})} \sum\limits_{t \in \mathcal{T}} \sum\limits_{i \in \mathcal{I}}  \alpha_i x_{ti} + \beta_i y_{ti} + \theta_t p_t + \eta_t s_t  \\
\st \quad  & z_i \in \{0,1\} \quad \forall \, i \in \mathcal{I} \\
& y_{ti} \in \{0,1\}  \quad  \forall \, i \in \mathcal{I}, \; t \in \mathcal{T}, \; \bm{d}^{[t]} \in \mathcal{D}^{[t]}, \, \bm{\hat{c}}^{[t-1]} \in \widehat{\mathcal{C}}^{[t-1]}(\bm{z}, \; \bm{y}^{[t-1]}) \\
& \hspace{-8pt} \left.
\begin{array}{l}
  s_t = s_{t-1} + \sum\limits_{i \in \mathcal{I}}x_{ti} + p_t - d_t \\[4pt]
  x_{ti} \leq c_i^{\max} - \hat{c}_{ti}  \quad  \forall \, i \in \mathcal{I} \\[4pt]
  x_{ti} \leq c_i^{\max} y_{ti} \quad  \forall \, i \in \mathcal{I} \\[4pt]
  x_{ti} \geq c_i^{\min} y_{ti}  \quad  \forall \, i \in \mathcal{I} \\[4pt]
  s_t \leq s^{\max} \\[4pt]
  s_t, \, p_t \in \mathbb{R}_+ \\[4pt]
  x_{ti} \in \mathbb{R}_+ \quad \forall \, i \in \mathcal{I}
\end{array}
\right\rbrace \quad \forall \, t \in \mathcal{T}, \; \bm{d}^{[t]} \in \mathcal{D}^{[t]}, \, \bm{\hat{c}}^{[t]} \in \widehat{\mathcal{C}}^{[t]}(\bm{z}, \; \bm{y}^{[t]})
\end{aligned}
\end{equation}
where $\bm{\alpha}$, $\bm{\beta}$, $\bm{\theta}$, and $\bm{\eta}$ are cost parameters, and $\bm{y}_t = \bm{y}_t(\bm{d}^{[t]}, \, \bm{\hat{c}}^{[t-1]})$, $\bm{x}_t = \bm{x}_t(\bm{d}^{[t]}, \, \bm{\hat{c}}^{[t]})$, $p_t = p_t(\bm{d}^{[t]}, \, \bm{\hat{c}}^{[t]})$, and $s_t = s_t(\bm{d}^{[t]}, \, \bm{\hat{c}}^{[t]})$. We consider the following uncertainty sets:
\begin{subequations}
\label{eqn:UncertaintySet-Multistage1}
\begin{align}
    & \mathcal{D}^{[t]} = \left \lbrace \bm{d}^{[t]} \in \mathbb{R}_+^{t}: \; \bm{d}^{\min} \leq \bm{d} \leq \bm{d}^{\max} \right \rbrace \\[4pt]
    & \widehat{\mathcal{C}}^{[t]}(\bm{z}, \; \bm{y}^{[t]}) = \left \lbrace 
    \bm{\hat{c}}^{[t]} \in \mathbb{R}_+^{tI}: 
    \begin{array}{l}
        % \hat{c}_{t'i} \geq \hat{c}_{t'i}^{\text{min-1}} z_i y_{t'i}(\bm{d}^{[t']}, \, \bm{\hat{c}}^{[t'-1]}) + \hat{c}_{t'i}^{\text{min-2}} (1-z_i) y_{t'i}(\bm{d}^{[t']}, \, \bm{\hat{c}}^{[t'-1]})  \quad \forall \, t' = 1, \ldots, t \\
        \hat{c}_{t'i} \leq \hat{c}_{t'i}^{\text{max-1}} z_i y_{t'i} + \hat{c}_{t'i}^{\text{max-2}} (1-z_i) y_{t'i}  \;\; \forall \,  t' = 1, \ldots, t, \, i \in \mathcal{I} \\[4pt]
        % \sum\limits_{t'=1}^t \sum\limits_{i \in \mathcal{I}} \hat{c}_{t'i} \leq \Gamma^t \sum\limits_{t'=1}^t \sum\limits_{i \in \mathcal{I}} \left( \hat{c}_{t'i}^{\text{max-1}} z_i y_{t'i} + \hat{c}_{t'i}^{\text{max-2}} (1-z_i) y_{t'i} \right )
        \sum\limits_{i \in \mathcal{I}} \hat{c}_{t'i} \leq \tau_{t'} \sum\limits_{i \in \mathcal{I}}  \hat{c}_{t'i}^{\text{max-1}} z_i y_{t'i} + \hat{c}_{t'i}^{\text{max-2}} (1-z_i) y_{t'i} \;\; \forall \,  t' = 1, \ldots, t
    \end{array}
    \right \rbrace,
    \label{eqn:DDUSinMultistageCase}
\end{align}
\end{subequations}
where $\mathcal{D}^{[t]}$ is a simple box uncertainty set, while $\widehat{\mathcal{C}}^{[t]}(\bm{z}, \; \bm{y}^{[t]})$ is a decision-dependent budget uncertainty set. The uncertain parameter $\hat{c}_{t'i}$ is nonzero only when $y_{t'i} = 1$. Furthermore, in the case of $y_{t'i} = 1$, the upper bound of $c_{t'i}$ is $c_{t'i}^{\text{max-1}}$ if $z_i = 1$ and $c_{t'i}^{\text{max-2}}$ otherwise. With $c_{t'i}^{\text{max-1}} < c_{t'i}^{\text{max-2}}$, an equipment upgrade increases the minimum capacity and reduces the level of uncertainty in the capacity. Note that there are bilinear terms involving two binary variables in \eqref{eqn:DDUSinMultistageCase}, i.e. $ z_i y_{t'i}$, which can be easily linearized.

In the following, we consider multiple instances of a production planning problem with three units. The maximum inventory $s^{\max}$ is set to 5; all the other data are provided in Appendix \ref{apx:Data}. We construct decision rules using one breakpoint for each uncertain parameter and such that recourse variables only depend on the uncertain parameters from the current time period.

\subsubsection{The Benefit of Discrete Recourse}

\begin{figure}[ht!]\centering
\includegraphics[width=5in]{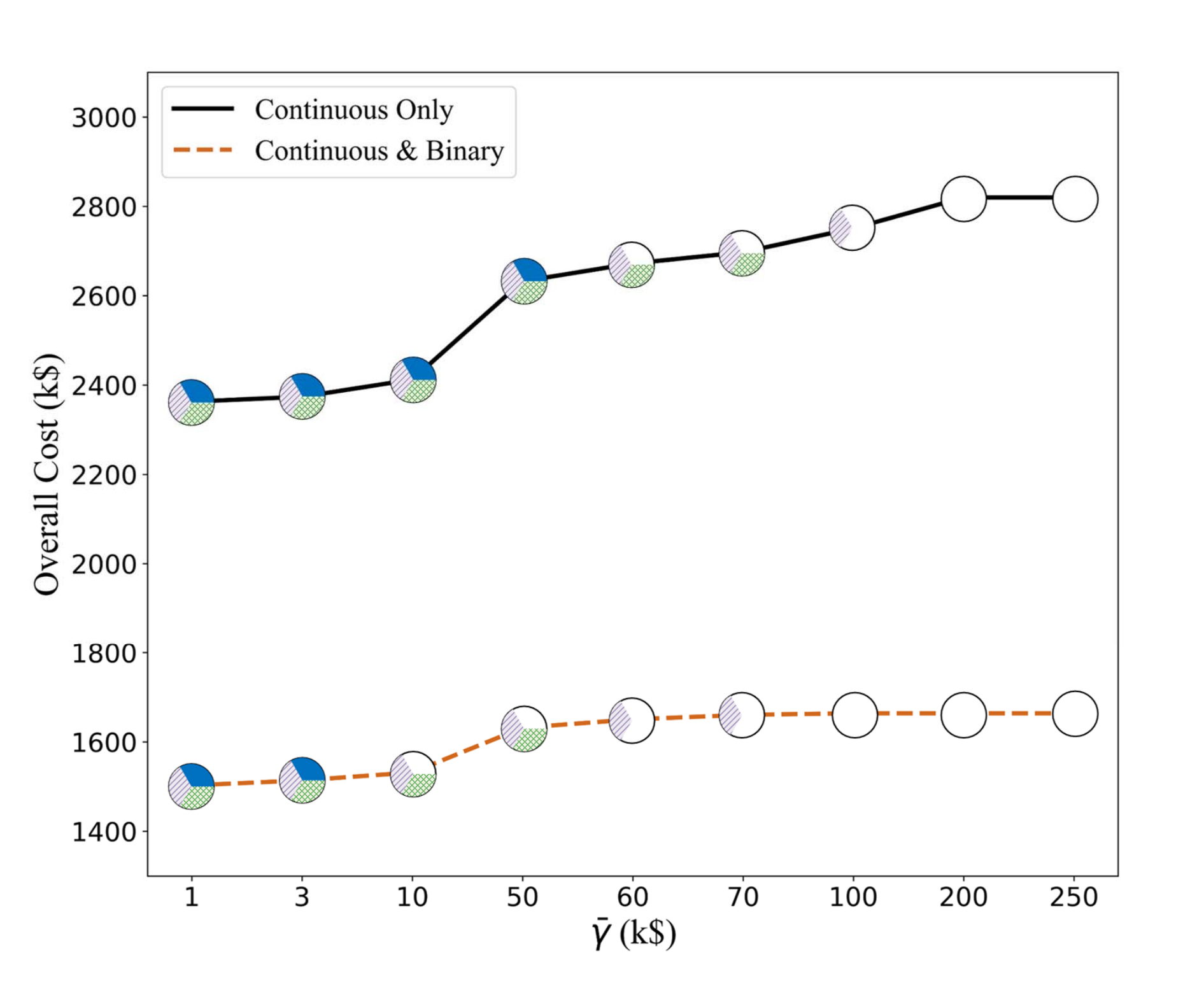}
\caption{Overall (worst-case) costs for different upgrade costs. Each instance is solved once with only continuous recourse and another time with both continuous and binary recourse. 
The purple stripes, green checks, and fill color blue indicate upgrades for Units 1, 2, and 3, respectively.
}
\label{fig:MultistageUpgradeCost}
\end{figure}

We first consider the case with two time periods and investigate how the solution depends on the equipment upgrade costs $\bm{\gamma}$. We set $\gamma_1 = \bar{\gamma}$, $\gamma_2 = 1.5 \bar{\gamma}$, and $\gamma_3 = 3 \bar{\gamma}$. The results from solving multiple instances with different $\bar{\gamma}$ are shown in Figure \ref{fig:MultistageUpgradeCost}. The pie chart associated with each instance indicates which units are being upgraded; colored fill means that we decide to upgrade the corresponding unit. One can see that as $\bar{\gamma}$ increases, it becomes less worthwhile to invest in equipment upgrades, up to a point where we leave all units unchanged. Figure \ref{fig:MultistageUpgradeCost} also shows the benefit of discrete recourse as we solve each instance once considering only continuous recourse and another time with both continuous and binary recourse. Comparing the optimal values, we see that a cost reduction of more than 35\,\% can be achieved if binary recourse is considered in addition to continuous recourse. Note that in both cases, continuous recourse variables can follow discontinuous piecewise linear decision rules.

We also solve the problem with different numbers of time periods. Table  \ref{tab:MultistageIntegerComparison} shows the results for $T$ equal to 2, 3, 4, and 5 with $\bar{\gamma} = 100$. Again, one can see that the cost is substantially reduced if in addition to continuous recourse, also binary recourse is considered. Moreover, the cost reduction increases with the number of time periods.

\begin{table}[h!]\centering
\setlength\tabcolsep{10pt}	
    \caption{Computational results for instances with different numbers of time periods. Each instance is solved once with only continuous recourse and another time with both continuous and binary recourse.}
    \begin{tabular}{llllll}
	\toprule
	& \# of time periods & 2   & 3   & 4   & 5 \\   
	\midrule
	\multirow{5}{*}{\makecell{\textbf{Continuous} \\ \textbf{recourse} \\ \textbf{only}}}  & Objective value (k\$)  & 2,751	&3,515	&4,032	&4,530 \\
    & \# of constraints  & 15,289	& 40,824	& 85,199	& 153,598
 \\
    & \# of continuous variables  & 4,647	& 11,793	& 23,883	& 42,213
 \\
    & \# of discrete variables  & 15 	& 21 & 27	& 33 
 \\
    & Computation time (s)  & 5	& 26	& 274	& 1,560
 \\
    \midrule
	\multirow{5}{*}{\makecell{\textbf{Continuous} \\ \textbf{\& binary} \\ \textbf{recourse}}}  & Objective value (k\$)  &1,664	& 2,112	& 2,434	& 2,745 \\
    & \# of constraints  & 69,805	& 203,175	& 443,663	& 822,373 \\
    & \# of continuous variables  & 18,777	& 53,571	& 115,671	& 212,853 \\
    & \# of discrete variables  & 51	& 75  & 99  & 123 \\
    & Computation time (s)  & 261	& 3,935	& 26,922	& 147,445 \\
    \bottomrule
	\end{tabular}
	\label{tab:MultistageIntegerComparison}
\end{table}

\subsubsection{Discussion on Computational Performance}

The computational results in Table \ref{tab:MultistageIntegerComparison} also show that the benefit of binary recourse comes at the cost of significantly greater computational complexity. For example, for $T=5$, the model with only continuous recourse solves in 1,560\,s, while the computation time for the model with both continuous and binary recourse is about two orders of magnitude longer. One reason for the higher computational complexity is obviously the increased model size. However, it turns out that the by far larger contributing factor is the ``looseness'' of the MILP formulation. Recall that the incorporation of binary recourse variables that affect the uncertainty set results in a formulation that involves bilinear terms. Each of these bilinear terms consists of a binary variable and a continuous variable representing a ``dual'' variable that stems from the reformulation. The linearization of these bilinear terms involves the lower bounds, which are zero, and the upper bounds of the dual variables. Generally, we cannot find tight upper bounds on these variables a priori such that very large values have to be chosen, which leads to a very weak LP relaxation of the MILP. This explanation is consistent with our observation, which is that the actual optimal solution is usually found fairly quickly but the lower bound only improves very slowly.

We further confirm our suspicion with a small experiment. After solving each of the instances with continuous and binary recourse to optimality, we update the bounds on the dual variables based on the optimal solution and re-solve the problem. Let $\Lambda$ be the value of a dual variable at the optimal solution, then the tightest update that we can apply is to set the upper bound $M$ to $\Lambda$. In addition, we consider two other update rules: $M = 2 \Lambda$ and $M = 2 \Lambda + 0.01$. The computation times for re-solving the four instances with the three different bound update rules are shown in Table \ref{tab:MultistageComputationTime}. In all cases, the computation times are drastically shorter than the ones required to solve the problems with the original bounds. Moreover, one can see that the solution times in case of the third update rule are significantly longer although the resulting bounds are only minimally larger. The reason is that at the optimal solution, most dual variables are zero; hence, the first two update rules fix all these variables to zero, which makes the problem considerably easier to solve. Note that these updated bounds do not result in a rigorous reformulation of the problem although the same optimal value is achieved. The sole purpose of this experiment is to demonstrate the impact of these bounds on the computation time.

\begin{table}[h!]\centering
\setlength\tabcolsep{10pt}	
    \caption{Computation times for re-solving instances with updated bounds in seconds.}
    \begin{tabular}{ccccc}
	\toprule
    \# of time periods  & 2   & 3   & 4   & 5 \\ 
    \midrule
    Update: $M = \Lambda$  & 1  & 4 & 8   & 12 \\
    Update: $M = 2 \Lambda$  & 2	& 4  	& 9  & 14 \\  
    Update: $M = 2 \Lambda + 0.01$  & 41	& 272	& 1,489   & 5,307 \\
    \bottomrule
	\end{tabular}
	\label{tab:MultistageComputationTime}
\end{table}

\section{Conclusions}
\label{sec:Conclusions}

In this work, we addressed multistage robust optimization with mixed-integer recourse and endogenous uncertainty, considering polyhedral uncertainty sets that are affected by binary variables. Applying a decision rule approach, which relies on the concept of lifted uncertainty sets, we derived tractable reformulations for the two- and multistage cases. The proposed framework has significant modeling flexibility as it can incorporate uncertainty sets affected by recourse decisions, binary recourse, and continuous recourse variables that follow discontinuous piecewise linear decision rules. The main advantage of appropriately modeling endogenous uncertainty and mixed-integer recourse is manifested in the significant reduction in solution conservatism, as demonstrated in our computational experiments.

Although the proposed reformulations enjoy the favorable tractability properties of robust optimization, the computational case studies also show that they tend to result in large and rather loose MILP formulations. Hence, future work will focus on the development of solution strategies that improve the computational performance.

\section*{Acknowledgments}
We gratefully acknowledge financial support from the National Key Research and Development Program of China (No. 2019YFB1705004), Science Fund for Creative Research Groups of NSFC (No. 61621002), and China Scholarship Council (CSC) (No. 201906320317).

\appendix
\section{Restricted Decision Rules in the Multistage Case}
\label{Appendix:DeltaTApproximation}

Consider decision rules for recourse variables in stage $t$ that only depend on uncertain parameters observed in the previous $\Delta t$ and current stages:
\begin{subequations}
\label{eqn:DeltaTApproximation-DecisionRule}
\begin{align}
    & \bm{x}_{t} = \bm{\overline{X}}_{11}^t \bm{\bar{\xi}}_{11} + \bm{\widehat{X}}_{11}^t \bm{\hat{\xi}}_{11} + \sum\limits_{t'= \max \{2, \, t - \Delta t \} }^t \sum\limits_{i=1}^{K_{t'}} \bm{\overline{X}}_{t'i}^t \bm{\bar{\xi}}_{t'i} + \bm{\widehat{X}}_{t'i}^t \bm{\hat{\xi}}_{t'i} \\
    & \bm{y}_t = \bm{\widehat{Y}}_{11}^t  \bm{\hat{\xi}}_{11} + \sum\limits_{t'= \max \{2, \, t - \Delta t \}}^t \sum\limits_{i=1}^{K_{t'}} \bm{\widehat{Y}}_{t'i}^t  \bm{\hat{\xi}}_{t'i},
\end{align}
\end{subequations}
where $\xi_{11}$ is still included in the decision rule as it accounts for the constant term. Following the same procedure presented in Subsection \ref{sec:MultiStageReformulation}, this results in the following reformulation of \eqref{eqn:MultiStageDeterministicConstraint} considering \eqref{eqn:DeltaTApproximation-DecisionRule}:
\begin{subequations}
\label{eqn:DeltaTApproximation-Reformulation}
\begin{align}
    % \minimize \quad & \bm{e}_1^{\top} \bm{x}_1 \\
    % \st \quad 
    % & \bm{A}^1 \bm{x}_1 + \bm{D}^1 \bm{y}_1 \leq \bm{b}^1 \\
    & \bm{f}_{t'i}^t = \bm{A}_{t'i}^t \bm{x}_1 + \bm{D}_{t'i}^t \bm{y}_1 - \bm{b}_{t'i}^t \quad \forall \, t= 2, \ldots, T, \, t' = 1, \ldots, t, \; i=1, \ldots, K_{t'} \\
    & \bm{\Phi}_t \bm{U}_1^t \bm{y}_1 + \sum\limits_{t'=1}^{t} \sum\limits_{i=1}^{K_{t'}} \bm{\delta}_{t'i}^t \leq \bm{0} \quad \forall \, t=2,\ldots, T \\
    & \left(\bm{f}_{11}^t  - \bm{\Phi}_{t} \bm{w}_{11}^{[t]} \right) v_{11} - \bm{\delta}_{11}^t \notag \\
    & + \sum\limits_{t'' = 2}^t \left( \bm{\widetilde{A}}_{t''}^t \bm{\overline{X}}_{11}^{t''} \bm{\bar{v}}_{11} + 
    % \sum\limits_{t'' = \max\{2, t' \}}^t 
    \left[ \bm{\widetilde{A}}_{t''}^t \bm{\widehat{X}}_{11}^{t''} + \left (\bm{\widetilde{D}}_{t''}^t  + \bm{\Phi}_t \bm{U}_{t''}^t \right) \left( \bm{\dot{\widehat{Y}}}_{11}^{t''} -\bm{\ddot{\widehat{Y}}}_{11}^{t''} \right)  \right] \bm{\hat{v}}_{11} \right) \leq \bm{0} \notag  \\
    & \qquad\qquad\qquad\qquad\qquad\qquad\qquad\qquad\qquad\qquad\qquad\qquad\qquad\qquad \forall \, t=2,\ldots,T \\
    & \left(\bm{f}_{t'i}^t  - \bm{\Phi}_{t} \bm{w}_{t'i}^{[t]} \right) v_{t'i} - \bm{\delta}_{t'i}^t \notag \\
    & + \sum\limits_{t'' = t'}^{\min \{t, \, t'+\Delta t \}} \left( \bm{\widetilde{A}}_{t''}^t \bm{\overline{X}}_{t'i}^{t''} \bm{\bar{v}}_{t'i} + 
    % \sum\limits_{t'' = \max\{2, t' \}}^t 
    \left[ \bm{\widetilde{A}}_{t''}^t \bm{\widehat{X}}_{t'i}^{t''} + \left (\bm{\widetilde{D}}_{t''}^t  + \bm{\Phi}_t \bm{U}_{t''}^t \right) \left( \bm{\dot{\widehat{Y}}}_{t'i}^{t''} -\bm{\ddot{\widehat{Y}}}_{t'i}^{t''} \right)  \right] \bm{\hat{v}}_{t'i} \right) \leq \bm{0}  \notag \\
    & \qquad\qquad\qquad \forall \, t=2,\ldots,T, \; t' = 2, \ldots, t, \; i = 1, \ldots, K_{t'}, \; \left( v_{t'i}, \, \bm{\bar{v}}_{t'i}, \, \bm{\hat{v}}_{t'i} \right) \in \widetilde{\mathcal{V}}_{t'i}.
\end{align}
\end{subequations}
The integrality constraints can be reformulated in the same fashion.

\section{Data for Case Studies}
\label{apx:Data}

\begin{table}[ht!]\centering
\setlength\tabcolsep{15pt}	
	\caption{Data for the two-stage design problem with three units.}
	\begin{tabular}{cccc}
	\toprule
	\textbf{Parameter} & \makecell{\textbf{Unit 1}} & \makecell{\textbf{Unit 2}} & \makecell{\textbf{Unit 3}} \\ \midrule
    $\alpha$ (k$\$$)     & 100      & 40     & 60 \\
    $\beta$ (k$\$$)       & 20       & 5      & 10     \\
    $\gamma$ (k$\$$)     & 5        & 2      & 4  \\
    $c^{\min}$        & 20           & 2          & 40  \\
    $c^{\max}$        & 145          & 65         & 140  \\
    $\hat{c}^{\max}$  & 35           & 20         & 5 \\
	\bottomrule
	\end{tabular}
	\label{appendix:data for twostage case with 3 alternative units}
\end{table}

\begin{table}[ht!]\centering
\setlength\tabcolsep{10pt}	
    \caption{Data for the two-stage design problem with eight units.}
	\begin{tabular}{ccccccccc}
	\toprule
	\textbf{Parameter} & \makecell{\textbf{Unit 1}} & \makecell{\textbf{Unit 2}} & \makecell{\textbf{Unit 3}} & \makecell{\textbf{Unit 4}} & \makecell{\textbf{Unit 5}} & \makecell{\textbf{Unit 6}} & \makecell{\textbf{Unit 7}} & \makecell{\textbf{Unit 8}} \\ \midrule
    $\alpha$ ($\$$)           & 75,365 & 61,420 & 98,153 & 66,932 & 81,824 & 62,627 & 83,175 & 66,110 \\
    $\beta$ ($\$$)           & 8,063  & 9,560  & 10,710  & 10,810  & 5,777  &  13,611  & 13,643  &  12,826     \\
    $\gamma$ ($\$$)          & 2,429 & 2,481  &  2,885   &  2,949   &  5,195   & 2,061   & 3,908   &   2,544  \\
    $c^{\min}$        & 21 & 13    & 20    & 28    & 2    & 11    & 30    & 21  \\
    $c^{\max}$        & 96   &  91   &  89   & 81   & 60   & 81   &  114   &  102  \\
    $\hat{c}^{\max}$  & 20.0  & 13.52  & 23.0  & 15.54  & 7.73  & 19.6  & 20.16  & 16.2 \\
	\bottomrule
	\end{tabular}
	\label{appendix:data for twostage case with 8 alternative units}
\end{table}

\begin{table}[ht!]\centering
\setlength\tabcolsep{10pt}	
    \caption{Cost and capacity data for the multistage production planning problem.}
    \begin{tabular}{cccc}
	\toprule
    \textbf{Parameter} & \makecell{\textbf{Unit 1}} & \makecell{\textbf{Unit 2}} & \makecell{\textbf{Unit 3}} \\ \midrule
    $\alpha$ (k$\$$)    &  2.0   & 3.0   & 5.5  \\
    $\beta$ (k$\$$)     & 20.0   & 40.0 & 80.0 \\
    $c^{\min}$      & 5.0     & 40.0   & 15.0  \\
    $c^{\max}$      & 50.0    & 100.0  & 90.0  \\
    \bottomrule
	\end{tabular}
	\label{appendix:design data for multistage production planning problem}
\end{table}

\begin{table}[ht!]\centering
\setlength\tabcolsep{10pt}	
    \caption{Data for time-varying uncertain parameters in the multistage production planning problem.}
    \begin{tabular}{ccccccc}
	\toprule
    \multicolumn{2}{c}{\makecell{\textbf{Parameter}}}   & \makecell{\textbf{Period 1}} & \makecell{\textbf{Period 2}} & \makecell{\textbf{Period 3}}   & \makecell{\textbf{Period 4}} & \makecell{\textbf{Period 5}}  \\ 
    \midrule
    \multicolumn{2}{c}{$\eta$ (k\$)}    & 3.5   & 3.5     & 3.5    & 3.5 & 3.5 \\
    \multicolumn{2}{c}{$\theta$ (k\$)}    & 15    &15     & 15    & 15     & 15 \\
    \midrule
    \multicolumn{2}{c}{$d^{\min}$}    & 35.0    & 54.0     & 33.0     & 27.0      & 25.0 \\
    \multicolumn{2}{c}{$d^{\max}$}    & 150.0   & 215.0    & 130.0    & 100.0     & 96.0 \\
    \midrule
     \multicolumn{2}{c}{$\tau $}    & 0.5   & 0.5   & 0.5   & 0.5   & 0.5   \\
    \midrule
    \multirow{3}{*}{$\hat{c}^{\text{max-1}}$}   & Unit 1     & 10.0   & 10.5  & 11.0  & 12.0  & 13.0  \\
    & Unit 2     & 5.0   & 6.0  & 7.0  & 7.5  & 8.0   \\
    & Unit 3     & 8.0   & 9.0  & 10.0  & 11.0  & 11.5   \\
    \midrule
    \multirow{3}{*}{$\hat{c}^{\text{max-2}}$}   & Unit 1     & 20.0   & 21.0  & 22.0  & 22.5  & 23.0   \\ 
    & Unit 2     & 15.0   & 16.0  & 17.0  & 18.0  & 18.5   \\ 
    & Unit 3     & 25.0   & 26.0  & 27.0  & 27.5  & 28.0   \\ 
    \bottomrule
	\end{tabular}
	\label{appendix:data for the time-varying uncertain parameters in the multistage production planning problem}
\end{table}

\bibliographystyle{newapa}
\bibliography{library}
%\bibliography{../../../../../References/BibTeX/library}

\end{document}